\documentclass[preprint,11pt,nopreprintline]{elsarticle}
\usepackage[margin=1in]{geometry}
\usepackage{amsmath,amssymb,amsthm}
\usepackage{indentfirst}
\usepackage[
    colorlinks=true,
    linkcolor=blue,
    citecolor=blue,
    urlcolor=blue
]{hyperref}

\hypersetup{pdfauthor={Qingying Deng and Chaoyang Zhang},
            pdftitle={Irreducibility of Full Multivariate Partial-Twuality Polynomials}}

\DeclareMathOperator{\rank}{rank}

\DeclareMathOperator{\im}{im}
\DeclareMathOperator{\corank}{corank}
\theoremstyle{plain}
\newtheorem{theorem}{Theorem}[section]
\newtheorem{lemma}[theorem]{Lemma}
\newtheorem{proposition}[theorem]{Proposition}
\newtheorem{corollary}[theorem]{Corollary}

\theoremstyle{definition}
\newtheorem{definition}[theorem]{Definition}
\newtheorem{example}[theorem]{Example}

\theoremstyle{remark}

\begin{document}
\begin{frontmatter}

\title{Irreducibility of Full Multivariate Partial-Twuality Polynomials}

\author[xtu]{Qingying Deng\corref{corresponding}}
\ead{qingying@xtu.edu.cn}
\author[xtu]{Chaoyang Zhang}
\ead{202521511267@smail.xtu.edu.cn}
\cortext[corresponding]{Corresponding author.}
\affiliation[xtu]{organization={Xiangtan University},
                 city={Xiangtan},
                 postcode={411105},
                 state={Hunan},
                 country={China}}

\begin{abstract}
For each nonempty finite simple graph, we study five full multivariate partial-twuality polynomials defined using matrix ranks and coranks over \(\mathbb F_2\). Over an arbitrary coefficient field $K$, we give a necessary and sufficient criterion for nontrivial factorization in the Laurent polynomial ring by additive separation of the exponent function across a nontrivial vertex partition. For $\bullet\in\{\tau,\delta\tau,\tau\delta,\delta\tau\delta\}$, the full multivariate partial-\(\langle \bullet \rangle\) polynomial is irreducible in \(K[z^{\pm1},\mathbf u,\mathbf v]\) if and only if the graph is connected. For \(\delta\), we characterize the connected graphs whose polynomials are reducible by a condition on the parities of vertex degrees and of the numbers of common neighbors, and determine their irreducible factors. For each operation, the corresponding polynomial is irreducible in \(K[z,\mathbf u,\mathbf v]\) if and only if the graph is connected and at least one state has zero \(z\)-exponent. In particular, the full multivariate partial-\(\langle\delta\rangle\) polynomial is irreducible in \(K[z,\mathbf u,\mathbf v]\) if and only if the graph is connected and bipartite.
\end{abstract}

\end{frontmatter}

\section{Introduction}

Geometric duality and Petrie duality generate a group of six operations isomorphic to $S_3$.
Wilson~\cite{Wilson1979} studied these operations for regular maps, and Lins~\cite{Lins1982} described the corresponding symmetries using graph-encoded maps. Writing $\delta$ and $\tau$ for geometric duality and Petrie duality, respectively, the five nonidentity
elements are $\delta,\tau,\delta\tau,\tau\delta,\delta\tau\delta$. We denote the fifth operation by \(\delta \tau \delta\); the relations in \(S_3\) give \(\delta \tau \delta =\tau \delta \tau\). Chmutov~\cite{Chmutov2009} introduced partial duality, which acts on a chosen subset of edges of a ribbon graph. Ellis-Monaghan and Moffatt~\cite{EllisMonaghanMoffatt2012Twisted} defined the ribbon-group action using partial duality and half-twists on individual edges. Gross, Mansour, and Tucker introduced the partial-dual Euler-genus polynomial in Part~I~\cite{GrossMansourTucker2020Distributions} and the remaining partial-twuality polynomials in Part~II~\cite{GrossMansourTucker2021Twuality}.

Deng, Jin, and Yan~\cite{DengJinYan2026} expressed the Euler genera arising from partial twuality of bouquets in terms of matrix ranks and coranks, and extended the five partial-twuality polynomials to square matrices. In this paper, we specialize their exponent functions to the adjacency matrix $M=M(G)$ of a finite simple graph $G$, regarded as a symmetric matrix over $\mathbb F_2$. Here, \(A\subseteq V(G)\) specifies a subset of the row and column indices of the adjacency matrix, whereas partial-twuality operations on ribbon graphs act on subsets of edges. The resulting multivariate polynomials are invariants under graph isomorphisms, with the vertex variables relabeled accordingly. In particular, the univariate $\langle\delta\rangle$-polynomial coincides with the partial-dual genus polynomial of simple graphs studied by Cheng~\cite{Cheng2025PartialDual}.

Let $G$ be a finite simple graph with nonempty vertex set $V$, and let $M$ be its adjacency matrix over $\mathbb F_2$. For each
\[
\bullet\in
\{\delta,\tau,\delta\tau,\tau\delta,\delta\tau\delta\},
\]
let $f_\bullet(A)$ denote the corresponding exponent function, computed from matrix ranks and coranks over $\mathbb F_2$. We consider the full multivariate refinement
\[
P_G^{\langle\bullet\rangle}(z;\mathbf u,\mathbf v)
=
\sum_{A\subseteq V}
\left(\prod_{x\in A}u_x\right)
\left(\prod_{x\in V\setminus A}v_x\right)
z^{f_\bullet(A)},
\]
where $z$ and the variables $\{u_x,v_x:x\in V\}$ are algebraically independent commuting indeterminates. The variable $u_x$ records that $x\in A$, whereas $v_x$ records that $x\notin A$. Setting $u_x=v_x=1$ for every $x\in V$ recovers the corresponding univariate matrix partial-twuality polynomial. Throughout the paper, the coefficient field $K$ is arbitrary and is independent of the field $\mathbb F_2$ over which the matrix ranks are computed.

Two members of this family are related to vertex-weighted graph polynomials. For Traldi's vertex-weighted graph bracket~\cite{TraldiWeightedBracket}, regard $G$ as unmarked, assign weights $v_x$ and $u_x$ to the unselected and selected states, respectively, and set the nullity variable to $z^{-1}$. Multiplication by $z^{|V|}$ then gives $P_G^{\langle\tau\rangle}$. For the weighted vertex-nullity interlace polynomial~\cite{TraldiWeightedInterlace}, add a loop at every vertex of $G$, assign weights $u_x$ and $v_x$ to the selected and unselected states, respectively, and set the factor raised to principal nullity to $z^{-1}$. Multiplication by $z^{\rank M}$ gives $P_G^{\langle\delta\tau\delta\rangle}$. Our contribution in these two cases concerns their factorization and irreducibility.

For each of the five polynomial families, disjoint
union gives a product of the corresponding polynomials.
Thus, disconnectedness produces reducibility, and we
ask when irreducibility characterizes connectedness.
A related phenomenon is known for interlace polynomials.
Chen, Jin, and Ma~\cite{ChenJinMa2026} proved that the
two-variable interlace polynomial of a nonempty
loopless graph is irreducible over $\mathbb C$ if
and only if the graph is connected.
They also obtained a connectivity characterization
for Courcelle's multivariate interlace polynomial
over $\mathbb C$; their Remark~4.3 observes that the multivariate proof works over an arbitrary field.
These results lead to the corresponding question for the five exponent functions considered here.
Irreducibility need not survive specialization, so the identities with weighted graph polynomials do not settle this question.
We determine the vertex partitions across which each exponent function is additively separable.

Our main tool is a general factorization criterion
that converts a nontrivial factorization of a
multiaffine subset-generating polynomial into an
additive separation of its exponent function across
a partition of the vertex set.
For
$\bullet\in\{\tau,\delta\tau,\tau\delta,\delta\tau\delta\}$,
we prove that $P_G^{\langle\bullet\rangle}$ is
irreducible in $K[z^{\pm1},\mathbf u,\mathbf v]$
if and only if $G$ is connected.

For $\bullet=\delta$, connectivity alone is insufficient,
and we give a complete characterization of the
connected graphs whose polynomials are reducible
in the Laurent polynomial ring.
Reducibility occurs precisely when complementing the
subgraph induced by the even-degree vertices produces
a disconnected graph and, for every pair of distinct
vertices, the number of common neighbors is odd if
and only if both vertices have even degree.
When these conditions hold, we give the complete
irreducible factorization in terms of the connected
components of the resulting graph.
In particular, for bipartite graphs, Laurent
irreducibility is equivalent to connectedness.

In the ordinary polynomial ring
$K[z,u_x,v_x:x\in V]$, we obtain a unified criterion
for all five operations:
$P_G^{\langle\bullet\rangle}$ is irreducible if and
only if $G$ is connected and
\[
\min_{A\subseteq V}f_\bullet(A)=0.
\]
For $\bullet=\delta$, the minimum is zero precisely
when $G$ is bipartite.
Thus $P_G^{\langle\delta\rangle}$ is irreducible in
the ordinary polynomial ring if and only if $G$
is connected and bipartite.
All our factorization and irreducibility results
hold over an arbitrary coefficient field $K$.

The remainder of the paper is organized as follows.
Section~2 recalls the partial-twuality polynomials
and their matrix formulations, and fixes the notation.
Section~3 establishes the general factorization
criteria and applies them to the five exponent
functions to obtain the irreducibility results
in the Laurent and ordinary polynomial rings.

\section{Preliminaries}
Deng, Jin, and Yan~\cite{DengJinYan2026} extended
partial-twuality polynomials to square matrices. In this paper, we consider the specialization of their definitions to adjacency matrices of finite simple graphs over \(\mathbb F_2\). We first fix the required notation.

Let \(G\) be a finite simple graph with vertex set \(V=V(G)\).
Its adjacency matrix over \(\mathbb F_2\) is
\[
M=(m_{xy})_{x,y\in V},
\]
where
\[
m_{xy}=
\begin{cases}
1, & xy\in E(G),\\
0, & \text{otherwise}.
\end{cases}
\]
For \(A\subseteq V\), write \(A^c=V\setminus A\), and let \(G[A]\) denote the subgraph of \(G\) induced by \(A\). Its adjacency matrix is the principal submatrix \(M[A]\) of \(M\) indexed by \(A\). For \(X,Y\subseteq V(G)\), let \(M[X,Y]\) denote the submatrix with rows indexed by \(X\) and columns indexed by \(Y\).
For a vector $b\in\mathbb F_2^V$ and a subset $A\subseteq V$, write $b[A]$
for its restriction to $A$. Let $\mathbf1_A,\mathbf0_A\in\mathbb F_2^A$
denote the all-ones and zero vectors, respectively.

Let \(I_A\) be the diagonal matrix indexed by \(V\) whose diagonal entry at \(x\) is \(1\) if \(x\in A\) and \(0\) otherwise. We write \(I_k\) for the identity matrix of order \(k\). Thus
\[
(M+I_A)[A]=M[A]+I_{|A|}.
\]
All matrix ranks, coranks, and determinants are taken over \(\mathbb F_2\). For a square matrix \(N\) of order \(k\), we write
\[
\corank N=k-\rank N.
\]
The symbols \(\ker N\) and \(\im N\) denote its kernel and column space, respectively.
For the empty matrix, we use the conventions
\[
\rank M[\varnothing]=\corank M[\varnothing]=0,
\qquad
\det M[\varnothing]=1.
\]
For $x\in V$, let
$N_G(x)=\{y\in V:xy\in E(G)\}$
denote its neighborhood. With this notation, the five matrix partial-twuality polynomials associated with a finite simple graph are defined as follows.

\begin{definition}[\cite{DengJinYan2026}]\label{def:graph-partial-twuality}
Let \(G\) be a finite simple graph with adjacency matrix \(M\) over \(\mathbb F_2\), indexed by \(V=V(G)\). For \(\bullet\in \{\delta,\tau,\delta\tau,\tau\delta,\delta\tau\delta\}\),
define
\[
P_G^{\langle\bullet\rangle}(z)
=
\sum_{A\subseteq V}
z^{r_{\langle\bullet\rangle}(M,A)},
\]
where, for each \(A\subseteq V\),
\[
\begin{aligned}
r_{\langle\delta\rangle}(M,A)
&=\rank M[A]+\rank M[A^c],\\
r_{\langle\tau\rangle}(M,A)
&=\rank(M+I_A),\\
r_{\langle\delta\tau\rangle}(M,A)
&=\rank\bigl((M+I_A)[A]\bigr)+\rank M[A^c],\\
r_{\langle\tau\delta\rangle}(M,A)
&=\rank(M+I_A)-\corank M[A],\\
r_{\langle\delta\tau\delta\rangle}(M,A)
&=\rank M-\corank\bigl((M+I_A)[A]\bigr).
\end{aligned}
\]
\end{definition}

By~\cite[Proposition~22]{DengJinYan2026}, all five exponent
functions in Definition~\ref{def:graph-partial-twuality} take
nonnegative integer values. In particular, the two expressions
involving subtraction are nonnegative because
\[\corank M[A]\leq\rank\bigl(M[A]+I_{|A|}\bigr)\leq\rank(M+I_A)
\]
and
\[\corank\bigl(M[A]+I_{|A|}\bigr)\leq\rank M[A]\leq\rank M.
\]

\section{Main Results}
Throughout this section, \(K\) denotes an arbitrary field. All graphs are finite and simple with at least one vertex. For a graph \(G\), we write \(V=V(G)\) and \(M=M(G)\).

\subsection{General factorization criteria}

Let
\[
\mathbf u={(u_x)}_{x\in V},
\qquad
\mathbf v={(v_x)}_{x\in V}
\]
be commuting indeterminates. For $A\subseteq V$, write
\[
u_A=\prod_{x\in A}u_x,
\qquad
v_A=\prod_{x\in A}v_x,
\]
with empty products equal to $1$.

\begin{definition}\label{def:full-multivariate}
For \(\bullet\in \{\delta,\tau,\delta\tau,\tau\delta,\delta\tau\delta\}\) and \(A\subseteq V\), set \(f_\bullet(A)=r_{\langle\bullet\rangle}(M,A)\), where \(r_{\langle\bullet\rangle}(M,A)\) is defined in Definition~\ref{def:graph-partial-twuality}. The full multivariate partial-\(\langle\bullet\rangle\) polynomial of \(G\) is defined by
\[
P_G^{\langle\bullet\rangle}(z;\mathbf u,\mathbf v)
=
\sum_{A\subseteq V}
u_A v_{V\setminus A}z^{f_\bullet(A)}.
\]
Setting \(u_x=v_x=1\) for every \(x\in V\) gives \(P_G^{\langle\bullet\rangle}(z;\mathbf 1,\mathbf 1)=P_G^{\langle\bullet\rangle}(z)\).
\end{definition}
When no confusion can arise, we suppress the variables and
write \(P_G^{\langle\bullet\rangle}\) for
\(P_G^{\langle\bullet\rangle}(z;\mathbf u,\mathbf v)\).

\begin{proposition}[\cite{DengJinYan2026}]\label{prop:block-diagonal-product}
Let \(G\) be a finite simple graph, and let \(V=V_1\sqcup V_2\) be a partition such that no edge joins \(V_1\) and \(V_2\). Then
\[
P_G^{\langle\bullet\rangle}(z)
=
P_{G[V_1]}^{\langle\bullet\rangle}(z)
P_{G[V_2]}^{\langle\bullet\rangle}(z)
\]
for every \(\bullet\in \{\delta,\tau,\delta\tau,\tau\delta,\delta\tau\delta\}\).
\end{proposition}

\begin{corollary}\label{cor:multivariate-disjoint-union}
Let \(G\) be a finite simple graph, and suppose that \(V=V_1\sqcup V_2\) is a partition such that no edge joins \(V_1\) and \(V_2\). Then
\[P_G^{\langle\bullet\rangle}=
P_{G[V_1]}^{\langle\bullet\rangle}
P_{G[V_2]}^{\langle\bullet\rangle}
\]
for every \(\bullet\in \{\delta,\tau,\delta\tau,\tau\delta,\delta\tau\delta\}\). In particular, if \(G=G_1\sqcup\cdots\sqcup G_k\) is the decomposition of \(G\) into its connected components,
then
\[P_G^{\langle\bullet\rangle}=\prod_{i=1}^{k}P_{G_i}^{\langle\bullet\rangle}.\]
\end{corollary}

\begin{proof}
Let \(M_i=M[V_i]\), and for \(A\subseteq V\), set
\(A_i=A\cap V_i\). Since no edge joins \(V_1\) and \(V_2\),
we have
\[
M=M_1\oplus M_2,
\qquad
I_A=I_{A_1}\oplus I_{A_2},
\]
where \(I_{A_i}\) is indexed by \(V_i\). It follows from the definitions of the five exponent functions
and the additivity of rank and corank under direct sums that
\[
r_{\langle\bullet\rangle}(M,A)
=
r_{\langle\bullet\rangle}(M_1,A_1)
+
r_{\langle\bullet\rangle}(M_2,A_2).
\]
Moreover,
\[
u_A v_{V\setminus A}
=
\bigl(u_{A_1}v_{V_1\setminus A_1}\bigr)
\bigl(u_{A_2}v_{V_2\setminus A_2}\bigr).
\]
Since every \(A\subseteq V\) is uniquely of the form
\(A=A_1\sqcup A_2\), summing over all
\(A_1\subseteq V_1\) and \(A_2\subseteq V_2\) gives the first
factorization. The factorization over the connected components
follows by repeated application of the two-part formula.
\end{proof}

The following lemma gives a general criterion for reducibility in a Laurent polynomial ring.

Let \(f:2^V\longrightarrow\mathbb Z_{\geq0}\), and define
\[\Phi_f=\sum_{A\subseteq V}u_Av_{A^c}z^{f(A)},\qquad A^c=V\setminus A.
\]

\begin{lemma}\label{lem:laurent-factorization-criterion}
The polynomial $\Phi_f$ is reducible in \(K[z^{\pm1},u_x,v_x:x\in V]\) if and only if there exists a nontrivial partition \(V=X\sqcup Y\) such that
\begin{equation}\label{eq:exponent-separation}
f(S\cup T)+f(\varnothing)=f(S)+f(T)
\qquad
(S\subseteq X,\ T\subseteq Y).
\end{equation}
\end{lemma}

\begin{proof}
Suppose first that there exists a nontrivial
partition $V=X\sqcup Y$ satisfying \eqref{eq:exponent-separation}. Since every $A\subseteq V$ has the unique representation \(A=(A\cap X)\cup(A\cap Y)\), we have
\[
\begin{aligned}
\Phi_f
&=
\sum_{S\subseteq X}\sum_{T\subseteq Y}
u_{S\cup T}v_{V\setminus(S\cup T)}
z^{f(S\cup T)}\\
&=
\sum_{S\subseteq X}\sum_{T\subseteq Y}
u_Sv_{X\setminus S}
u_Tv_{Y\setminus T}
z^{f(S)+f(T)-f(\varnothing)}\\
&=
z^{-f(\varnothing)}
\left(
\sum_{S\subseteq X}
u_Sv_{X\setminus S}z^{f(S)}
\right)
\left(
\sum_{T\subseteq Y}
u_Tv_{Y\setminus T}z^{f(T)}
\right).
\end{aligned}
\]
The units of
$K[z^{\pm1},u_x,v_x:x\in V]$ are precisely the elements
$cz^m$, where $c\in K^\times$ and $m\in\mathbb Z$.
Since $X,Y\neq\varnothing$, neither sum is a unit.
The factor $z^{-f(\varnothing)}$ is a unit and may be
absorbed into either sum. Hence $\Phi_f$ is reducible.

Conversely, suppose that $\Phi_f$ is reducible in
$K[z^{\pm1},u_x,v_x:x\in V]$.
Since $z^{-f(\varnothing)}$ is a unit, we may write
\[
\Phi_f=z^{-f(\varnothing)}FH,
\]
where $F,H\in K[z^{\pm1},u_x,v_x:x\in V]$
are nonunits.
Both factors are nonzero since $\Phi_f\ne0$.

Fix $x\in V$. Every monomial of $\Phi_f$ contains exactly
one of $u_x$ and $v_x$. Hence $\Phi_f$ is homogeneous of
total degree $1$ in the pair $(u_x,v_x)$. Regarding all
the remaining variables as coefficients, the additivity
of total degree over an integral domain gives
\[
\deg_{u_x,v_x}F+\deg_{u_x,v_x}H=1.
\]
Since these degrees are nonnegative integers, exactly one
is $1$ and the other is $0$.

Suppose first that
\[
\deg_{u_x,v_x}F=1,
\qquad
\deg_{u_x,v_x}H=0.
\]
Then $H$ is independent of $u_x$ and $v_x$.
Setting $u_x=v_x=0$ in
$\Phi_f=z^{-f(\varnothing)}FH$ gives
\[
0
=
\left.\Phi_f\right|_{u_x=v_x=0}
=
z^{-f(\varnothing)}
\left.F\right|_{u_x=v_x=0}H.
\]
Since $H\neq0$ and the Laurent polynomial ring in the remaining
variables is an integral domain, it follows that
\[
\left.F\right|_{u_x=v_x=0}=0.
\]
Thus $F$ has no term of degree $0$ in $u_x,v_x$. Since its
total degree in these variables is $1$, $F$ is homogeneous
of degree $1$ in $(u_x,v_x)$.

The case
\[
\deg_{u_x,v_x}F=0,
\qquad
\deg_{u_x,v_x}H=1
\]
is analogous. Moreover, both $u_x$ and $v_x$ occur in the
factor of degree $1$, since the other factor is independent
of these variables, while $\Phi_f$ contains nonzero terms
involving each of them. Therefore, for every $x\in V$, the
pair $(u_x,v_x)$ occurs in exactly one of $F$ and $H$.

Thus $F$ and $H$ induce a partition $V=X\sqcup Y$,
where $(u_x,v_x)$ occurs in $F$ for $x\in X$
and in $H$ for $x\in Y$.
Both $X$ and $Y$ are nonempty.
If $X=\varnothing$, then $F\in K[z^{\pm1}]$.
Taking the coefficient of $v_V$ in
$\Phi_f=z^{-f(\varnothing)}FH$ shows that
$F$ divides $z^{2f(\varnothing)}$.
Since this is a unit, $F$ is also a unit,
a contradiction.
The case $Y=\varnothing$ is analogous.

Comparing the coefficients of $v_V$ shows that the
product of the coefficients of $v_X$ in $F$ and
$v_Y$ in $H$ is $z^{2f(\varnothing)}$.
Hence both coefficients are units in $K[z^{\pm1}]$.
After multiplying $F$ and $H$ by mutually inverse units,
we may assume that both coefficients are
$z^{f(\varnothing)}$.

For every $S\subseteq X$, the coefficient of
$u_Sv_{X\setminus S}v_Y$ in $\Phi_f$ equals the
coefficient of $u_Sv_{X\setminus S}$ in $F$:
the factor $z^{-f(\varnothing)}$ cancels the
coefficient $z^{f(\varnothing)}$ of $v_Y$ in $H$.
Therefore,
\[
F=
\sum_{S\subseteq X}
u_Sv_{X\setminus S}z^{f(S)}.
\]
Similarly, comparing the coefficients of
$v_Xu_Tv_{Y\setminus T}$ for every $T\subseteq Y$ gives
\[
H=
\sum_{T\subseteq Y}
u_Tv_{Y\setminus T}z^{f(T)}.
\]

Now fix $S\subseteq X$ and $T\subseteq Y$.
The coefficient of
$u_{S\cup T}v_{V\setminus(S\cup T)}$
in $z^{-f(\varnothing)}FH$ is
$z^{f(S)+f(T)-f(\varnothing)}$,
whereas its coefficient in $\Phi_f$ is
$z^{f(S\cup T)}$.
Therefore,
\[
z^{f(S\cup T)}
=
z^{f(S)+f(T)-f(\varnothing)}.
\]
Equality of these Laurent monomials gives
\[
f(S\cup T)+f(\varnothing)=f(S)+f(T).
\]

\end{proof}

\begin{corollary}\label{cor:ordinary-factorization-criterion}
The polynomial $\Phi_f$ is reducible in \(K[z,u_x,v_x:x\in V]\) if and only if at least one of the following conditions holds:
\begin{enumerate}
\item
\[
\min_{A\subseteq V}f(A)>0;
\]
\item there exists a nontrivial partition \(V=X\sqcup Y\) such that
\[
f(S\cup T)+f(\varnothing)=f(S)+f(T)
\qquad
(S\subseteq X,\ T\subseteq Y).
\]
\end{enumerate}
\end{corollary}

\begin{proof}
Suppose first that
\[
\min_{A\subseteq V}f(A)>0.
\]
Then
\[
\Phi_f
=
z\left(
\sum_{A\subseteq V}
u_Av_{V\setminus A}z^{f(A)-1}
\right).
\]
All exponents in the sum are nonnegative, so the sum
belongs to $K[z,u_x,v_x:x\in V]$.
Since $V\ne\varnothing$, it has total degree $|V|>0$
in the variables $u_x,v_x$ and is therefore a nonunit.
Since $z$ is also a nonunit in this ring, the displayed
expression is a nontrivial factorization of $\Phi_f$.

It remains to consider the case
\[
\min_{A\subseteq V}f(A)=0,
\]
in which $z\nmid\Phi_f$.
The ring $K[z,u_x,v_x:x\in V]$ is a unique
factorization domain, and the Laurent polynomial ring
$K[z^{\pm1},u_x,v_x:x\in V]$ is obtained from it
by inverting $z$.
Under this localization, the nonzero constant multiples
of $z$ become units, while every other irreducible
polynomial remains an irreducible nonunit.
Since $z\nmid\Phi_f$, none of the irreducible factors
of $\Phi_f$ is a nonzero constant multiple of $z$.
Consequently, $\Phi_f$ is reducible in the ordinary
polynomial ring if and only if it is reducible
in the Laurent polynomial ring.

By Lemma~\ref{lem:laurent-factorization-criterion},
the latter holds if and only if there exists a
nontrivial partition $V=X\sqcup Y$ satisfying
\eqref{eq:exponent-separation}.
Thus, when $\min_{A\subseteq V}f(A)=0$, reducibility
is equivalent to condition~(2).
Together with the first case, this proves the stated
criterion.
\end{proof}

\subsection{Separation properties and Laurent irreducibility}

\subsubsection{The partial-\texorpdfstring{\(\langle\tau\rangle\)}{<tau>} polynomial}

\begin{lemma}[\cite{DengJinYan2026}]\label{lem:nonsingular-diagonal-perturbation}
Let \(M\) be an \(n\times n\) matrix over \(\mathbb F_2\), with rows and columns indexed by a finite set \(V\). Then there exists \(A\subseteq V\) such that \(M+I_A\) is nonsingular.
\end{lemma}

For $A\subseteq V$, let
\[
f_\tau(A)=\rank(M+I_A).
\]

\begin{lemma}\label{lem:tau-separation}
Let \(G\) be a finite simple graph with adjacency matrix \(M\) over \(\mathbb{F}_2\), and let \(V=X\sqcup Y\) be a nontrivial partition. Then
\begin{equation}\label{eq:tau-separation}
f_\tau(S\cup T)+f_\tau(\varnothing)=f_\tau(S)+f_\tau(T)
\qquad(S\subseteq X,\ T\subseteq Y)
\end{equation}
if and only if no edge of \(G\) joins \(X\) and \(Y\).
\end{lemma}

\begin{proof}
Set \(n=|V|\). If no edge joins $X$ and $Y$, then, after ordering the
vertices with those in $X$ first, \(M=M[X]\oplus M[Y]\). Consequently, for \(S\subseteq X\) and \(T\subseteq Y\),
\[
M+I_{S\cup T}
=
\bigl(M[X]+I_S\bigr)\oplus
\bigl(M[Y]+I_T\bigr),
\]
where \(I_S\) and \(I_T\) are indexed by \(X\) and \(Y\),
respectively. By the additivity of rank under direct sums,
\[
\begin{aligned}
f_\tau(S\cup T)+f_\tau(\varnothing)
&=\rank(M[X]+I_S)+\rank(M[Y]+I_T)\\
&\quad+\rank M[X]+\rank M[Y],\\
f_\tau(S)+f_\tau(T)
&=\rank(M[X]+I_S)+\rank M[Y]\\
&\quad+\rank M[X]+\rank(M[Y]+I_T).
\end{aligned}
\]
The two expressions are equal, proving~\eqref{eq:tau-separation}.

Conversely, suppose that~\eqref{eq:tau-separation} holds. By Lemma~\ref{lem:nonsingular-diagonal-perturbation}, there exists \(A\subseteq V\) such that \(B=M+I_A\) is nonsingular, and put \(C=B^{-1}=(c_{xy})_{x,y\in V}\). Since \(M\) is symmetric, both \(B\) and \(C\) are symmetric. Set \(S=A\cap X\), and \(T=A\cap Y\).

For $x\in X$ and $y\in Y$, condition~\eqref{eq:tau-separation}
gives
\[
\begin{aligned}
f_\tau(A)
&=f_\tau(S)+f_\tau(T)-f_\tau(\varnothing),\\
f_\tau(A\triangle\{x\})
&=f_\tau(S\triangle\{x\})+f_\tau(T)-f_\tau(\varnothing),\\
f_\tau(A\triangle\{y\})
&=f_\tau(S)+f_\tau(T\triangle\{y\})-f_\tau(\varnothing),\\
f_\tau(A\triangle\{x,y\})
&=f_\tau(S\triangle\{x\})
 +f_\tau(T\triangle\{y\})
 -f_\tau(\varnothing).
\end{aligned}
\]
Hence
\[
f_\tau(A\triangle\{x,y\})+f_\tau(A)
=
f_\tau(A\triangle\{x\})+f_\tau(A\triangle\{y\}),
\]
where $\triangle$ denotes symmetric difference. Since \(B=M+I_A\) is nonsingular, this gives
\[
\rank(B+I_{\{x,y\}})+n
=
\rank(B+I_{\{x\}})
+
\rank(B+I_{\{y\}}).
\]
All three ranks are at most \(n\). Hence \(B+I_{\{x,y\}}\) is nonsingular if and only if both \(B+I_{\{x\}}\) and \(B+I_{\{y\}}\) are nonsingular. Over \(\mathbb{F}_2\), this gives
\begin{equation}\label{eq:tau-determinant-identity}
\det\bigl(B+I_{\{x,y\}}\bigr)=\det\bigl(B+I_{\{x\}}\bigr)\det\bigl(B+I_{\{y\}}\bigr).
\end{equation}

Since \(\det C=1\), multiplying these matrices on the left by \(C\) and expanding along the unchanged columns of the identity matrix yields
\[
\det\bigl(B+I_{\{x\}}\bigr)=1+c_{xx},
\qquad
\det\bigl(B+I_{\{y\}}\bigr)=1+c_{yy}.
\]
Similarly, using the symmetry of \(C\),
\[
\begin{aligned}
\det\bigl(B+I_{\{x,y\}}\bigr)
&=
\det
\begin{pmatrix}
1+c_{xx}&c_{xy}\\
c_{xy}&1+c_{yy}
\end{pmatrix}
&=(1+c_{xx})(1+c_{yy})+c_{xy}^{\,2}.
\end{aligned}
\]
Substituting into~\eqref{eq:tau-determinant-identity} gives
\(c_{xy}^{\,2}=0\), and hence \(c_{xy}=0\).

Since \(x\in X\) and \(y\in Y\) were arbitrary, \(C\) is block
diagonal with respect to \(X\sqcup Y\), and so is \(B=C^{-1}\).
Since \(M\) and \(B\) differ only on the diagonal,
\(M[X,Y]=0\). Thus no edge joins \(X\) and \(Y\).
\end{proof}

\begin{theorem}\label{thm:tau-irreducibility-laurent}
Let $G$ be a finite simple graph with at least one vertex. Then $P_G^ {\langle\tau\rangle}(z;\mathbf u,\mathbf v)$ is irreducible in $ K[z^{\pm1},u_x,v_x:x\in V]$ if and only if $G$ is connected.
\end{theorem}

\begin{proof}
By Lemma~\ref{lem:laurent-factorization-criterion},
$P_G^{\langle\tau\rangle}(z;\mathbf u,\mathbf v)$ is reducible
if and only if there is a nontrivial partition $V=X\sqcup Y$
such that
\[
f_{\tau}(S\cup T)+f_{\tau}(\varnothing)=f_{\tau}(S)+f_{\tau}(T)
\]
for every $S\subseteq X$ and $T\subseteq Y$.
By Lemma~\ref{lem:tau-separation}, such a partition exists if and
only if $G$ is disconnected. The result follows.
\end{proof}

\subsubsection{The partial-\texorpdfstring{\(\langle\delta\tau\rangle\)}{<delta tau>} polynomial}

For $A\subseteq V$, let
\[
f_{\delta\tau}(A)=\rank\bigl((M+I_A)[A]\bigr)+\rank M[A^c].
\]

\begin{theorem}\label{thm:delta-tau-irreducibility-laurent}
Let \(G\) be a finite simple graph with at least one vertex. Then \(P_G^{\langle\delta\tau\rangle}(z;\mathbf u,\mathbf v)\) is irreducible in \(K[z^{\pm1},u_x,v_x:x\in V]\) if and only if \(G\) is connected.
\end{theorem}

\begin{proof}
If $G$ is disconnected, Corollary~\ref{cor:multivariate-disjoint-union}
expresses $P_G^{\langle\delta\tau\rangle}$ as a product over its
connected components. Each component is nonempty, so its factor has
positive degree in the vertex variables and is a nonunit. Hence the
polynomial is reducible.

Conversely, suppose that \(G\) is connected and that \(P^{\langle \delta \tau \rangle}_G\) is reducible. By Lemma~\ref{lem:laurent-factorization-criterion}, there exists a nontrivial partition \(V=X\sqcup Y\) such that
  \[f_{\delta \tau}(S\cup T)+f_{\delta \tau}(\varnothing)=f_{\delta \tau}(S)+f_{\delta \tau}(T) \qquad (S\subseteq X,T\subseteq Y).\]
  Since \(G\) is connected, some edge \(xy\) joins \(X\) and \(Y\). Taking \(S=\{x\}\) and \(T=\{y\}\), we obtain

  \begin{equation}\label{eq:delta-tau-parity}
      f_{\delta \tau}(\{x,y\})+f_{\delta \tau}(\varnothing)=f_{\delta \tau}(\{x\})+f_{\delta \tau}(\{y\}).
  \end{equation}

  Since \(M\) is the adjacency matrix of a simple graph, every principal submatrix of \(M\) is symmetric with zero diagonal and hence is alternating over \(\mathbb F_2\). An alternating bilinear form is an orthogonal sum of its radical and
two-dimensional nondegenerate alternating forms; see
\cite[Lemma~2.3]{Aliabadi2026}. The alternating condition excludes
one-dimensional blocks with nonzero diagonal entries. Each nondegenerate
block therefore has rank two, so every such submatrix has even rank. Consequently, \(f_{\delta \tau}(\varnothing)=\rank M\) is even, whereas
  \[f_{\delta \tau}(\{x\})=1+\rank M[V\setminus \{x\}],\quad f_{\delta \tau}(\{y\})=1+\rank M[V\setminus \{y\}]\]
  are odd. Moreover, since \(xy\) is an edge,
  \[M[\{x,y\}]+I_2=\begin{pmatrix}
    1&1\\1&1
  \end{pmatrix}\]
  has rank \(1\). Thus \(f_{\delta \tau}(\{x,y\})=1+\rank M[V\setminus \{x,y\}]\) is odd. The left-hand side of~\eqref{eq:delta-tau-parity} is therefore odd, while the right-hand side is even, a contradiction. Hence \(P^{\langle \delta \tau \rangle}_G\) is irreducible.

\end{proof}

\subsubsection{The partial-\texorpdfstring{\(\langle\delta\tau\delta\rangle\)}{<delta tau delta>} polynomial}

For $A\subseteq V$, let
\[
f_{\delta\tau\delta}(A)
=
\rank M-\corank\bigl((M+I_A)[A]\bigr).
\]

\begin{theorem}\label{thm:delta-tau-delta-laurent-irreducibility}
Let \(G\) be a finite simple graph with at least one vertex. Then \(P_G^{\langle\delta\tau\delta\rangle}(z;\mathbf u,\mathbf v)\) is irreducible in \(K[z^{\pm1},u_x,v_x:x\in V]\) if and only if \(G\) is connected.
\end{theorem}

\begin{proof}
If $G$ is disconnected, Corollary~\ref{cor:multivariate-disjoint-union}
expresses $P_G^{\langle\delta\tau\delta\rangle}$ as a product over its
connected components. Each component is nonempty, so its factor has
positive degree in the vertex variables and is a nonunit. Hence the
polynomial is reducible.

Conversely, suppose that \(G\) is connected and that
\(P_G^{\langle\delta\tau\delta\rangle}\) is reducible. By Lemma~\ref{lem:laurent-factorization-criterion}, there exists a nontrivial partition \(V=X\sqcup Y\) such that
\begin{equation}\label{eq:tau-delta-tau-separation}
f_{\delta \tau \delta}(S\cup T)+f_{\delta \tau \delta}(\varnothing)=f_{\delta \tau \delta}(S)+f_{\delta \tau \delta}(T)
\qquad(S\subseteq X,\ T\subseteq Y).
\end{equation}

Since \(M\) has zero diagonal,
\[
M[\{x\}]+I_1=[1]
\qquad(x\in V).
\]
Hence
\begin{equation}\label{eq:tau-delta-tau-singleton}
f_{\delta \tau \delta}(\varnothing)=f_{\delta \tau \delta}(\{x\})=\rank M
\qquad(x\in V).
\end{equation}

Fix \(x\in X\) and \(y\in Y\). Writing \(m_{xy}\) for the
corresponding entry of \(M\), we have
\[
M[\{x,y\}]+I_2
=
\begin{pmatrix}
1&m_{xy}\\
m_{xy}&1
\end{pmatrix}.
\]
This matrix has corank \(0\) if \(m_{xy}=0\), and corank \(1\)
if \(m_{xy}=1\). Therefore
\begin{equation}\label{eq:tau-delta-tau-pair}
f_{\delta \tau \delta}(\{x,y\})=\rank M-m_{xy},
\end{equation}
where \(m_{xy}\) is regarded as the integer \(0\) or \(1\).

Applying~\eqref{eq:tau-delta-tau-separation} with
\(S=\{x\}\) and \(T=\{y\}\), and using
equations~\eqref{eq:tau-delta-tau-singleton} and~\eqref{eq:tau-delta-tau-pair}
gives
\[
\bigl(\rank M-m_{xy}\bigr)
+\rank M=2\rank M.
\]
Thus \(m_{xy}=0\). Since \(x\in X\) and \(y\in Y\) were
arbitrary, no edge joins \(X\) and \(Y\), contradicting the
connectedness of \(G\). Therefore
\(P_G^{\langle \delta\tau\delta \rangle}\) is irreducible.
\end{proof}

\subsubsection{The partial-\texorpdfstring{\(\langle\tau\delta\rangle\)}{<tau delta>} polynomial}

For $A\subseteq V$, let
\[
f_{\tau\delta}(A)
=
\rank(M+I_A)-\corank M[A].
\]

\begin{lemma}[\cite{Aliabadi2026}] \label{lem:rank-one-update}
Let $H$ be a symmetric matrix over $\mathbb F_2$, and let
$b$ be a column vector of the appropriate size. Then
\[
\rank(H+bb^{\mathsf T})
=
\begin{cases}
\rank H+1,
& b\notin\im H,\\
\rank H,
& Hq=b,\ b^{\mathsf T}q=0,\\
\rank H-1,
& Hq=b,\ b^{\mathsf T}q=1.
\end{cases}
\]
In the last two cases, $q$ is any solution of $Hq=b$, and
the value of $b^{\mathsf T}q$ is independent of the chosen
solution.
\end{lemma}

\begin{lemma}\label{lem:alternating-rank}
Let $M$ be a symmetric matrix with zero diagonal over
$\mathbb F_2$, indexed by a finite set $V$.
For $A\subseteq V$, write
$r=\rank M$,
$\rho(A)=\rank(M+I_A)$.
For $x\in V$, let $e_x\in\mathbb F_2^V$ denote the standard
basis vector corresponding to $x$.

Then
\begin{equation}\label{eq:rho-single}
\rho(\{x\})=
\begin{cases}
r, & e_x\in\im M,\\
r+1, & e_x\notin\im M
\end{cases}
\qquad(x\in V).
\end{equation}
Moreover, for distinct $x,y\in V$ with
$e_x\notin\im M$,
\begin{equation}\label{eq:rho-pair}
\rho(\{x,y\})=r
\quad\Longleftrightarrow\quad
e_x+e_y\in\im M.
\end{equation}
\end{lemma}

\begin{proof}
Fix an ordering of $V$.
For every $q\in\mathbb F_2^V$, the symmetry and zero diagonal
of $M$ imply that
\[
q^{\mathsf T}Mq
=
\sum_{x\in V}m_{xx}q_x^2
+
\sum_{x<y}(m_{xy}+m_{yx})q_xq_y
=0.
\]

We have $I_{\{x\}}=e_xe_x^{\mathsf T}$.
Apply Lemma~\ref{lem:rank-one-update} with $H=M$ and $b=e_x$.
If $e_x\notin\im M$, then
\[
\rho(\{x\})=\rank(M+I_{\{x\}})=r+1.
\]
If $e_x\in\im M$, choose $q$ such that $Mq=e_x$.
Then $e_x^{\mathsf T}q=q^{\mathsf T}Mq=0$.
Lemma~\ref{lem:rank-one-update} therefore gives
\[
\rho(\{x\})=\rank(M+I_{\{x\}})=r.
\]
This proves \eqref{eq:rho-single}.

Now let $x\ne y$ and assume that
$e_x\notin\im M$.
Put $H=M+I_{\{x\}}$.
By \eqref{eq:rho-single}, $\rank H=r+1$.

Suppose first that $\rho(\{x,y\})=r$.
The identity
\[
M+I_{\{x,y\}}=H+e_ye_y^{\mathsf T}
\]
shows that the update at $y$ decreases the rank of $H$ by one.
Lemma~\ref{lem:rank-one-update} provides $q$ such that \(Hq=e_y\), \(e_y^{\mathsf T}q=q_y=1\). Multiplying $Hq=e_y$ on the left by $q^{\mathsf T}$, we obtain
\[
q_y=q^{\mathsf T}e_y=q^{\mathsf T}Hq=q^{\mathsf T}Mq+q_x^2=q_x.
\]
Here $q^{\mathsf T}Mq=0$ and $q_x^2=q_x$ in $\mathbb F_2$.
Hence $q_x=q_y=1$, and \(Mq=e_y+q_xe_x=e_x+e_y\). Therefore, $e_x+e_y\in\im M$.

Conversely, suppose that $e_x+e_y\in\im M$,
and choose $q$ such that $Mq=e_x+e_y$.
Then
\[
0=q^{\mathsf T}Mq=q^{\mathsf T}(e_x+e_y)=q_x+q_y,
\]
so $q_x=q_y$.
The symmetry of $M$ gives
$\im M=(\ker M)^\perp$.
As $e_x\notin\im M$, there exists
$h\in\ker M$ such that $e_x^{\mathsf T}h\ne0$.
The only nonzero element of $\mathbb F_2$ is $1$, so
$h_x=e_x^{\mathsf T}h=1$.
Furthermore,
\[
0=h^{\mathsf T}(e_x+e_y)=h_x+h_y,
\]
since $e_x+e_y\in\im M$ and $h\in\ker M$.
It follows that $h_y=1$.

The equality $q_x=q_y$ leaves two possibilities: \(q_x=q_y=1\) or \(q_x=q_y=0\). In the second case, put $q'=q+h$.
Then
\[
Mq'=M(q+h)=Mq+Mh=e_x+e_y,
\]
and
\[
q'_x=q_x+h_x=1,\qquad
q'_y=q_y+h_y=1.
\]
Thus, in either case, we can choose a solution $q$ satisfying
\[
Mq=e_x+e_y,\qquad q_x=q_y=1.
\]
It follows that
\[
Hq=(M+I_{\{x\}})q
=Mq+q_xe_x
=e_x+e_y+e_x
=e_y,
\qquad e_y^{\mathsf T}q=q_y=1.
\]
Applying Lemma~\ref{lem:rank-one-update} to $H$ and $e_y$ gives
\[
\rank(H+e_ye_y^{\mathsf T})
=\rank H-1=r.
\]
Finally, $H+e_ye_y^{\mathsf T}=M+I_{\{x,y\}}$, and hence
$\rho(\{x,y\})=r$.
This proves \eqref{eq:rho-pair}.
\end{proof}

\begin{lemma}\label{lem:alternating-deletion}
With the notation of Lemma~\ref{lem:alternating-rank}, we have
\begin{equation}\label{eq:deletion-rank}
\rank M[V\setminus\{x\}]
=
\begin{cases}
r-2, & e_x\in\im M,\\
r, & e_x\notin\im M
\end{cases}
\qquad(x\in V).
\end{equation}
If \(e_x\notin\im M\), inserting a zero in coordinate \(x\) induces a linear isomorphism from \(\ker M[V\setminus\{x\}]\) onto \(H_x:=\{h\in\ker M:h_x=0\}\).
\end{lemma}

\begin{proof}
Fix $x\in V$.
Since $M$ is symmetric, $\im M=(\ker M)^\perp$.

Suppose first that $e_x\in\im M$.
Then every vector in $\ker M$ has zero $x$-coordinate.
Restriction to $V\setminus\{x\}$ therefore defines an
injective linear map from $\ker M$ into
$\ker M[V\setminus\{x\}]$.
Hence
\[
\dim\ker M[V\setminus\{x\}]
\ge \dim\ker M
=|V|-r.
\]
Since $M[V\setminus\{x\}]$ has order $|V|-1$, we obtain
\[
\rank M[V\setminus\{x\}]
=|V|-1-\dim\ker M[V\setminus\{x\}]
\le r-1.
\]
Deleting one row and the corresponding column decreases
rank by at most two, so
\[
r-2\le\rank M[V\setminus\{x\}]\le r-1.
\]
Both $M$ and $M[V\setminus\{x\}]$ are symmetric with
zero diagonal over $\mathbb F_2$, and hence have even rank.
Therefore,
\[
\rank M[V\setminus\{x\}]=r-2.
\]

Now suppose that $e_x\notin\im M$.
Choose $h\in\ker M$ with $h_x=1$.
For $q\in\ker M[V\setminus\{x\}]$, insert a zero in
coordinate $x$, leaving all other coordinates unchanged,
and continue to denote the resulting vector by $q$.
Since $(Mq)_v=0$ for every $v\ne x$ and $h_x=1$,
\[
h^{\mathsf T}Mq
=\sum_{v\in V}h_v(Mq)_v
=h_x(Mq)_x
=(Mq)_x.
\]
On the other hand, symmetry and $Mh=0$ give
\[
h^{\mathsf T}Mq=(Mh)^{\mathsf T}q=0.
\]
Thus $Mq=0$, so the extended vector belongs to $\ker M$
and has $x$-coordinate zero.

Conversely, restriction to \(V\setminus\{x\}\) maps \(H_x\) into \(\ker M[V\setminus\{x\}]\) and is the inverse of the zero extension map. This proves the asserted linear isomorphism.

Since the coordinate functional at $x$ is nonzero on
$\ker M$, the subspace \(H_x\) has dimension $\dim\ker M-1$.
The isomorphism therefore gives
\[
\begin{aligned}
\corank M[V\setminus\{x\}]
&=\dim \ker M[V\setminus \{x\}]\\
&=\dim H_x=\dim \ker M-1=\corank M-1,
\end{aligned}
\]
and hence
\[
\rank M[V\setminus\{x\}]=r.
\]
\end{proof}

\begin{lemma}\label{lem:tau-delta-cross-edge-characterization}
Let $G$ be a finite simple graph with adjacency matrix
$M$ over $\mathbb F_2$, and let $V=X\sqcup Y$ be a
nontrivial partition. For $A\subseteq V$, define
\(f_{\tau\delta}(A)=\rank(M+I_A)-\corank M[A]\).
Suppose that
\begin{equation}\label{eq:tau-delta-cross-edge-separation}
f_{\tau\delta}(S\cup T)+f_{\tau\delta}(\varnothing)
=
f_{\tau\delta}(S)+f_{\tau\delta}(T)
\qquad(S\subseteq X,\ T\subseteq Y).
\end{equation}
Then, for every $x\in X$ and $y\in Y$,
\[
xy\in E(G)
\quad\Longleftrightarrow\quad
e_x\notin\im M
\text{ and }e_x+e_y\in\im M.
\]
\end{lemma}

\begin{proof}
Write $r=\rank M$, and fix $x\in X$ and $y\in Y$.
Let $m_{xy}\in\{0,1\}$ be the corresponding entry of $M$,
viewed as an integer.

For \(A\subseteq V\), write \(\rho(A)=\rank(M+I_A)\).
For \(x\in X\) and \(y\in Y\), we have
\[
\corank M[\{x\}]=\corank M[\{y\}]=1,
\qquad
\corank M[\{x,y\}]=2-2m_{xy}.
\]
Since \(f_{\tau\delta}(\varnothing)=r\), taking
\(S=\{x\}\) and \(T=\{y\}\) in
\eqref{eq:tau-delta-cross-edge-separation} gives
\begin{equation}\label{eq:general-tau-delta-cross-rank}
\rho(\{x,y\})-\rho(\{x\})-\rho(\{y\})+r
=
-2m_{xy}.
\end{equation}

Suppose first that \(xy\in E(G)\). Then
\[
\rho(\{x,y\})-\rho(\{x\})
=
\rho(\{y\})-r-2.
\]
Since
\[
M+I_{\{x,y\}}
=
\bigl(M+I_{\{x\}}\bigr)+e_ye_y^{\mathsf T},
\]
applying Lemma~\ref{lem:rank-one-update} with
\(H=M+I_{\{x\}}\) and \(b=e_y\) gives
\[
\rho(\{x,y\})-\rho(\{x\})\ge -1.
\]
Thus \(\rho(\{y\})-r-2\ge -1\), and hence
\(\rho(\{y\})\ge r+1\).
The single-vertex rank formula in
Lemma~\ref{lem:alternating-rank} now gives
\(\rho(\{y\})=r+1\).
Interchanging \(x\) and \(y\), we also obtain
\(\rho(\{x\})=r+1\).
Substituting these values into
\eqref{eq:general-tau-delta-cross-rank} yields
\(\rho(\{x,y\})=r\).
The single-vertex rank formula in
Lemma~\ref{lem:alternating-rank} implies that
\(e_x\notin\im M\), and its two-vertex rank formula then gives
\(e_x+e_y\in\im M\).

Conversely, suppose that \(e_x\notin\im M\),
\(e_x+e_y\in\im M\).
We also have \(e_y\notin\im M\): otherwise,
\(e_x=(e_x+e_y)+e_y\in\im M\), contrary to the assumption.
The rank formulas in Lemma~\ref{lem:alternating-rank} give
\[
\rho(\{x\})=\rho(\{y\})=r+1,
\qquad
\rho(\{x,y\})=r.
\]
Substitution into
\eqref{eq:general-tau-delta-cross-rank} yields \(-2=-2m_{xy}\).
Thus \(m_{xy}=1\), and hence \(xy\in E(G)\).
\end{proof}

\begin{lemma}\label{lem:tau-delta-separation}
Let \(G\) be a finite simple graph with adjacency matrix \(M\) over \(\mathbb F_2\), and let \(V(G)=X\sqcup Y\) be a nontrivial partition. Then
\begin{equation}\label{eq:tau-delta-separation}
f_{\tau\delta}(S\cup T)+f_{\tau\delta}(\varnothing)
=
f_{\tau\delta}(S)+f_{\tau\delta}(T)
\qquad
(S\subseteq X,\ T\subseteq Y)
\end{equation}
holds if and only if no edge of \(G\) joins \(X\) and \(Y\).
\end{lemma}

\begin{proof}
Suppose first that no edge joins \(X\) and \(Y\). After
ordering the vertices with those in \(X\) first, we have
\(M=M[X]\oplus M[Y]\). The additivity of rank and corank under
direct sums gives, for every \(A\subseteq V\),
\[
f_{\tau\delta}(A)+f_{\tau\delta}(\varnothing)
=
f_{\tau\delta}(A\cap X)+f_{\tau\delta}(A\cap Y).
\]
Taking \(A=S\cup T\), where \(S\subseteq X\) and
\(T\subseteq Y\), yields~\eqref{eq:tau-delta-separation}.

Conversely, assume that~\eqref{eq:tau-delta-separation} holds
and that an edge joins \(X\) and \(Y\). The same identity holds on the connected component containing this
edge. Indeed, for subsets of that component, the block-diagonal form
of \(M\) shows that every other component contributes the same fixed
rank to \(\rank(M+I_A)\) and does not enter \(M[A]\); hence these
contributions cancel from the two sides of
\eqref{eq:tau-delta-separation}. Replacing \(G\) by this component,
and \(X,Y\) by their intersections with its vertex set, we may therefore assume that \(G\) is connected.

Write
\[
n=|V|,\qquad
r=\rank M,\qquad
\rho(A)=\rank(M+I_A).
\]
Applying Lemma~\ref{lem:tau-delta-cross-edge-characterization}, we obtain
\begin{equation}\label{eq:cross-edge-characterization}
xy\in E(G)
\quad\Longleftrightarrow\quad
e_x\notin\im M
\ \text{and}\
e_x+e_y\in\im M
\qquad
(x\in X,\ y\in Y).
\end{equation}

We first show that an edge joining \(X\) and \(Y\) forces the kernel
to be spanned by the all-ones vector and every vertex of \(X\) to be
adjacent to every vertex of \(Y\). Fix such an edge \(xy\),
with \(x\in X\) and \(y\in Y\), and set
\[
W=\{v\in V:e_v+e_x\in\im M\}.
\]
Since \(e_x+e_x=0\in\im M\), we have \(x\in W\).
Equation~\eqref{eq:cross-edge-characterization} gives \(y\in W\).
Thus \(W\) meets both \(X\) and \(Y\).

Every \(v\in W\) satisfies \(e_v\notin\im M\). Indeed, if
\(e_v\in\im M\), then
\[
e_x=(e_v+e_x)+e_v\in\im M,
\]
contrary to~\eqref{eq:cross-edge-characterization} applied to
the edge \(xy\).

We claim that no edge joins \(W\) and \(V\setminus W\).
Assume that \(ij\) is such an edge, with \(i\in W\) and
\(j\notin W\). If \(i\) and \(j\) lie in different parts of
\(X\sqcup Y\), then~\eqref{eq:cross-edge-characterization}
gives \(e_i+e_j\in\im M\). Consequently,
\[
e_j+e_x
=
(e_i+e_j)+(e_i+e_x)
\in\im M,
\]
contradicting \(j\notin W\). Hence \(i\) and \(j\) lie in the same part.

The case \(i,j\in Y\) is analogous, using \(x\in W\cap X\).
We may therefore assume that \(i,j\in X\). As \(i,y\in W\),
\[
e_i+e_y
=
(e_i+e_x)+(e_y+e_x)
\in\im M.
\]
On the other hand, \(e_i+e_j\notin\im M\); otherwise,
\[
e_j+e_x
=
(e_i+e_j)+(e_i+e_x)
\in\im M,
\]
again contradicting \(j\notin W\).

We already know that \(e_i\notin\im M\).
Lemma~\ref{lem:alternating-rank} therefore gives
\[
\rho(\{i\})=r+1,\qquad
\rho(\{i,y\})=r,\qquad
\rho(\{i,j\})\ne r.
\]
A rank-one update changes rank by at most one. It follows that
\[
\rho(\{i,j\})\ge \rho(\{i\})-1=r.
\]
The inequality is strict because \(\rho(\{i,j\})\ne r\), so
\(\rho(\{i,j\})\ge r+1\). Similarly,
\(\rho(\{i,y,j\})\le\rho(\{i,y\})+1=r+1\).

The edge \(ij\) implies that
\(\rank M[\{i,j\}]=2\), and hence
\(\rank M[\{i,y,j\}]\ge2\). The latter matrix is symmetric
with zero diagonal over \(\mathbb F_2\), so its rank is even.
Therefore,
\[
\rank M[\{i,y,j\}]=2,
\qquad
\corank M[\{i,y,j\}]=1.
\]

We now have
\[
f_{\tau\delta}(\{i,j\})
=
\rho(\{i,j\})
\ge r+1,
\]
whereas
\[
f_{\tau\delta}(\{i,y,j\})
=
\rho(\{i,y,j\})-1
\le r.
\]
Furthermore, \(y\in W\) gives \(e_y\notin\im M\), so
\[
f_{\tau\delta}(\{y\})
=
\rho(\{y\})-1
=
r
=
f_{\tau\delta}(\varnothing).
\]
Applying~\eqref{eq:tau-delta-separation} with
\(S=\{i,j\}\) and \(T=\{y\}\) gives
\[
f_{\tau\delta}(\{i,y,j\})
=
f_{\tau\delta}(\{i,j\}),
\]
a contradiction. Thus no edge joins \(W\) and
\(V\setminus W\).

The graph \(G\) is connected and \(W\ne\varnothing\), so
\(W=V\). By the definition of \(W\) and the preceding
observation, every \(v\in V\) satisfies \(e_v+e_x\in\im M\) and \(e_v\notin\im M\). Moreover, for any \(v,j\in V\),
\[
e_v+e_j
=
(e_v+e_x)+(e_j+e_x)
\in\im M.
\]

Let \(\mathbf 1\) denote the all-ones vector, and take
\(h\in\ker M\). The symmetry of \(M\) gives
\(\im M=(\ker M)^\perp\). Therefore \(h_v=h_x\) for every \(v\in V\), and hence \(h=h_x\mathbf 1\). As \(h\in\ker M\) was arbitrary, we obtain \(\ker M\subseteq\langle\mathbf 1\rangle\). On the other hand, \(e_x\notin\im M\) implies that \(M\) is singular, so \(\ker M\ne\{0\}\). We conclude that
\begin{equation}\label{eq:kernel-structure}
\ker M=\langle\mathbf 1\rangle,
\qquad
r=n-1.
\end{equation}

Together with~\eqref{eq:cross-edge-characterization}, the
preceding observations show that every vertex of \(X\) is
adjacent to every vertex of \(Y\). Hence
\begin{equation}\label{eq:complete-cross-block}
M[X,Y]=J_{X,Y},
\end{equation}
where \(J_{X,Y}\) denotes the all-ones matrix.

We now show that the structure described by
\eqref{eq:kernel-structure} and
\eqref{eq:complete-cross-block} is incompatible
with~\eqref{eq:tau-delta-separation}. The matrix \(M\) is
symmetric with zero diagonal over \(\mathbb F_2\), so \(r\)
is even. As \(r=n-1\), the integer \(n\) is odd. By
interchanging \(X\) with \(Y\), and simultaneously \(x\)
with \(y\), if necessary, we may assume that \(|X|\) is even
and \(|Y|\) is odd. As the partition is nontrivial,
\(|X|\ge 2\).

Set \(B=M+I_{\{y\}}\). We have \(e_y\notin\im M\), so
\(\rank B=r+1=n\) by~\eqref{eq:rho-single}. Thus \(B\)
is nonsingular. Write
\[
C=B^{-1}=(c_{uv})_{u,v\in V}.
\]
The matrix \(C\) is symmetric since \(B\) is symmetric.

For \(A\subseteq V\), order the vertices with those in
\(A\) first. Then
\[
(B+I_A)C
=
I+I_AC
=
\begin{pmatrix}
I_{|A|}+C[A] & C[A,V\setminus A]\\
0            & I_{n-|A|}
\end{pmatrix}.
\]
Eliminating the upper-right block gives
\begin{equation}\label{eq:inverse-rank-identity}
\rank(B+I_A)
=
n-|A|
+
\rank\bigl(I_{|A|}+C[A]\bigr).
\end{equation}

Let \(u\in X\). By the preceding conclusions, \(e_u\notin\im M\) and \(e_u+e_y\in\im M\). Hence \eqref{eq:rho-pair} gives
\(\rho(\{u,y\})=r=n-1\). As \(B+I_{\{u\}}=M+I_{\{u,y\}}\), applying
\eqref{eq:inverse-rank-identity} with \(A=\{u\}\) gives
\(\rank\bigl([1+c_{uu}]\bigr)=0\). Hence \(c_{uu}=1\).

Now let \(u,v\in X\) be distinct. Applying
\eqref{eq:inverse-rank-identity} with \(A=\{u,v\}\) gives
\[
\rho(\{u,v,y\})
=n-2+
\rank
\begin{pmatrix}
0      & c_{uv}\\
c_{uv} & 0
\end{pmatrix}
=
n-2+2c_{uv},
\]
where \(c_{uv}\in\{0,1\}\) is viewed as an integer in the final
expression. By~\eqref{eq:complete-cross-block}, both \(uy\)
and \(vy\) are edges. The symmetric zero-diagonal matrix
\(M[\{u,v,y\}]\) therefore has rank \(2\) and corank \(1\).
It follows that
\[
f_{\tau\delta}(\{u,v,y\})
=
n-3+2c_{uv}.
\]

We also have \(e_u\notin\im M\) and
\(e_u+e_v\in\im M\). Equation~\eqref{eq:rho-pair} gives
\(\rho(\{u,v\})=n-1\). Writing \(m_{uv}\) for the corresponding
entry of \(M\), we obtain
\[
f_{\tau\delta}(\{u,v\})
=
(n-1)-(2-2m_{uv})
=
n-3+2m_{uv}.
\]
Moreover,
\(f_{\tau\delta}(\{y\})=f_{\tau\delta}(\varnothing)=r\).
Applying~\eqref{eq:tau-delta-separation} with
\(S=\{u,v\}\) and \(T=\{y\}\) makes the preceding two
expressions equal. Thus \(c_{uv}=m_{uv}\). Together with
\(c_{uu}=1\), this proves
\begin{equation}\label{eq:inverse-principal-block}
C[X]=M[X]+I_{|X|}.
\end{equation}

Finally,~\eqref{eq:kernel-structure} gives
\(M\mathbf1=0\). Taking the \(X\)-coordinates with respect to
\(V=X\sqcup Y\), we obtain
\[
M[X]\mathbf1_X+M[X,Y]\mathbf1_Y=0.
\]
Equation~\eqref{eq:complete-cross-block} and the oddness of
\(|Y|\) give
\[
M[X,Y]\mathbf1_Y
=
J_{X,Y}\mathbf1_Y
=
\mathbf1_X.
\]
Hence \(M[X]\mathbf1_X=\mathbf1_X\).

With respect to the ordering \(V(G)=X\sqcup Y\), set
\[
w=
\begin{pmatrix}
\mathbf1_X\\
\mathbf0_Y
\end{pmatrix}.
\]
The evenness of \(|X|\), together with
\eqref{eq:complete-cross-block}, gives
\[
Mw
=
\begin{pmatrix}
M[X]\mathbf1_X\\
J_{Y,X}\mathbf1_X
\end{pmatrix}
=
\begin{pmatrix}
\mathbf1_X\\
\mathbf0_Y
\end{pmatrix}
=
w.
\]
As \(y\in Y\), we have \(I_{\{y\}}w=0\), and hence \(Bw=w\).
Multiplying by \(C=B^{-1}\) gives \(Cw=w\). In particular,
\(C[X]\mathbf1_X=\mathbf1_X\). On the other hand,
\eqref{eq:inverse-principal-block} gives
\[
C[X]\mathbf1_X
=
\bigl(M[X]+I_{|X|}\bigr)\mathbf1_X
=
\mathbf1_X+\mathbf1_X
=
\mathbf0_X.
\]
This is impossible because \(X\ne\varnothing\). Therefore no
edge joins \(X\) and \(Y\).
\end{proof}

\begin{theorem}\label{thm:tau-delta-irreducibility}
Let \(G\) be a finite simple graph with at least one vertex. Then \(P_G^{\langle\tau\delta\rangle}(z;\mathbf u,\mathbf v)\) is irreducible in \(K[z^{\pm1},u_x,v_x:x\in V]\) if and only if \(G\) is connected.
\end{theorem}

\begin{proof}
If $G$ is disconnected, Corollary~\ref{cor:multivariate-disjoint-union}
expresses $P_G^{\langle\tau\delta\rangle}$ as a product over its
connected components. Each component is nonempty, so its factor has
positive degree in the vertex variables and is a nonunit. Hence the
polynomial is reducible.

Conversely, suppose that \(G\) is connected and that
\(P_G^{\langle \tau\delta \rangle}\) is reducible. By Lemma~\ref{lem:laurent-factorization-criterion}, there exists a
nontrivial partition \(V=X\sqcup Y\) such that
\[
f_{\tau\delta}(S\cup T)+f_{\tau\delta}(\varnothing)
=
f_{\tau\delta}(S)+f_{\tau\delta}(T)
\qquad
(S\subseteq X,\ T\subseteq Y).
\]
Lemma~\ref{lem:tau-delta-separation} then implies that no edge joins \(X\) and \(Y\), contrary
to the connectedness of \(G\). Hence \(P_G^{\langle \tau\delta \rangle}\) is
irreducible.
\end{proof}

\subsubsection{The partial-\texorpdfstring{\(\langle\delta\rangle\)}{<delta>} polynomial}

For $A\subseteq V$, let
\[
f_\delta(A)=\rank M[A]+\rank M[A^c].
\]

\begin{lemma}\label{lem:cross-edge-characterization}
Let \(G\) be a finite simple graph with adjacency matrix \(M\) over
\(\mathbb F_2\), and let \(V=X\sqcup Y\) be a nontrivial partition.
Suppose that, for every \(x\in X\) and \(y\in Y\),
\begin{equation}\label{eq:cross-edge-separation}
f_\delta(\{x,y\})+f_\delta(\varnothing)
=
f_\delta(\{x\})+f_\delta(\{y\}).
\end{equation}
Then, for every \(x\in X\) and \(y\in Y\),
\[
xy\in E(G)
\quad\Longleftrightarrow\quad
e_x\notin\operatorname{im}M
\ \text{and}\
e_x+e_y\in\operatorname{im}M.
\]
\end{lemma}

\begin{proof}
Write $r=\rank M$, and fix $x\in X$ and $y\in Y$.
Let $m_{xy}\in\{0,1\}$ be the corresponding entry of $M$,
viewed as an integer. Since \(M\) has zero diagonal,
\[
\rank M[\{x\}]=\rank M[\{y\}]=0,
\qquad
\rank M[\{x,y\}]=2m_{xy}.
\]
By \eqref{eq:cross-edge-separation} and \(f_\delta(\varnothing)=r\), we obtain
\begin{equation}\label{eq:delta-cross-deletion}
2m_{xy}+\rank M[V\setminus\{x,y\}]+r
=\rank M[V\setminus\{x\}]+\rank M[V\setminus\{y\}].
\end{equation}

Suppose first that \(xy\in E(G)\), so \(m_{xy}=1\).
Deleting one row and the corresponding column decreases rank
by at most two. Hence
\[
\rank M[V\setminus\{x,y\}]
\ge
\rank M[V\setminus\{x\}]-2.
\]
Combining this inequality with
\eqref{eq:delta-cross-deletion} gives
\[
\rank M[V\setminus\{y\}]\ge r.
\]
Since a principal submatrix has rank at most \(r\), equality holds.
Interchanging \(x\) and \(y\) gives
\[
\rank M[V\setminus\{x\}]
=
\rank M[V\setminus\{y\}]
=
r.
\]
By Lemma~\ref{lem:alternating-deletion},
\[
e_x\notin\im M,
\qquad
e_y\notin\im M.
\]
Equation~\eqref{eq:delta-cross-deletion} now yields
\[
\rank M[V\setminus\{x,y\}]=r-2.
\]

Deleting the row and column indexed by \(y\) from
\(M[V\setminus\{x\}]\) gives \(M[V\setminus\{x,y\}]\).
Since
\[
\rank M[V\setminus\{x\}]=r,
\qquad
\rank M[V\setminus\{x,y\}]=r-2,
\]
this deletion decreases rank by two.
Applying Lemma~\ref{lem:alternating-deletion}
to $M[V\setminus\{x\}]$ at $y$, we obtain
\(e_y\in\im M[V\setminus\{x\}]\).
For any $h\in\ker M$ with $h_x=0$, deleting its zero
$x$-coordinate gives a vector in $\ker M[V\setminus\{x\}]$,
whose $y$-coordinate is $h_y$.
Since $M[V\setminus\{x\}]$ is symmetric and
$e_y\in\im M[V\setminus\{x\}]$, we obtain $h_y=0$.

Interchanging \(x\) and \(y\), the same argument shows that
\(h_y=0\) implies \(h_x=0\) for every \(h\in\ker M\).
Thus \(h_x\) and \(h_y\) are either both zero or both nonzero.
Since they lie in \(\mathbb F_2\), we obtain \(h_x=h_y\)
for every \(h\in\ker M\).
For every \(h\in\ker M\), we have
\[
(e_x+e_y)^{\mathsf T}h=h_x+h_y=0,
\]
where the last equality holds over \(\mathbb F_2\).
Thus \(e_x+e_y\in(\ker M)^\perp=\im M\),
since \(M\) is symmetric.

Conversely, suppose that \(e_x\notin\im M\),
\(e_x+e_y\in\im M\).
Then \(e_y\notin\im M\); otherwise,
\(e_x=(e_x+e_y)+e_y \in \im M\).
Lemma~\ref{lem:alternating-deletion} gives
\[
\rank M[V\setminus\{x\}]=\rank M[V\setminus\{y\}]=r.
\]
Moreover, \(e_x+e_y\in\im M=(\ker M)^\perp\) implies that
\(h_x=h_y\) for every \(h\in\ker M\).
By Lemma~\ref{lem:alternating-deletion}, every vector in
\(\ker M[V\setminus\{x\}]\) therefore has zero \(y\)-coordinate.
Consequently, \(e_y\in\im M[V\setminus\{x\}]\).
Applying Lemma~\ref{lem:alternating-deletion} once more gives
\[
\rank M[V\setminus\{x,y\}]=r-2.
\]
Substituting these ranks into \eqref{eq:delta-cross-deletion}
yields \(2m_{xy}+(r-2)+r=2r\), so \(m_{xy}=1\),
and hence \(xy\in E(G)\).
\end{proof}

\begin{lemma}\label{lem:delta-rank-one-cut}
Let $G$ be a connected finite simple graph with adjacency matrix $M$ over $\mathbb F_2$, and let $V=X\sqcup Y$ be a nontrivial partition. Suppose that
\begin{equation}\label{eq:delta-structure-additivity}
f_\delta(S\cup T)+f_\delta(\varnothing)
=
f_\delta(S)+f_\delta(T)
\qquad(S\subseteq X,\ T\subseteq Y).
\end{equation}
Then $\corank M=1$. If $h$ is the nonzero vector in $\ker M$, then \(M[X,Y]=h[X]h[Y]^{\mathsf T}\).
\end{lemma}

\begin{proof}
Fix a vertex $y\in Y$ incident with an edge joining $X$ and $Y$. By Lemma~\ref{lem:cross-edge-characterization}, $e_y\notin\im M$. Lemma~\ref{lem:alternating-deletion} therefore gives
\[
f_\delta(\{y\})=f_\delta(\varnothing).
\]
Taking \(T=\{y\}\) in \eqref{eq:delta-structure-additivity}, and using
\(f_\delta(\{y\})=f_\delta(\varnothing)\), we obtain \(f_\delta(S\cup\{y\})=f_\delta(S)\quad(S\subseteq X)\).

Expanding this equality gives
\begin{equation}\label{eq:delta-root-increment}
\begin{aligned}
\rank M[S\cup\{y\}]-\rank M[S]=\rank M[V\setminus S]-\rank M[V\setminus(S\cup\{y\})].
\end{aligned}
\end{equation}

We claim that every vector in \(\ker M[V\setminus\{y\}]\) vanishes on each component of \(G[X]\) that contains a neighbor of \(y\).

Suppose otherwise. Choose \(h\in\ker M[V\setminus\{y\}]\) that is nonzero on a component of \(G[X]\) containing a neighbor of \(y\). Choose a shortest path in \(G[X\cup\{y\}]\) from \(y\) to a vertex where \(h\) is nonzero. Denote its terminal vertex by \(x\). Let \(S\) be the set of internal vertices of this path.
Then \(h_v=0\) for every \(v\in S\); otherwise, ending
the path at such a vertex would give a shorter path
to a nonzero coordinate of \(h\).

By Lemma~\ref{lem:alternating-deletion}, inserting a zero at \(y\) extends \(h\) to a vector in \(\ker M\). Its restrictions to the two matrices on the right-hand side of \eqref{eq:delta-root-increment} are kernel vectors with \(x\)-coordinate one. By symmetry and Lemma~\ref{lem:alternating-deletion}, deleting \(x\) preserves both ranks. Thus the right-hand side is unchanged when \(S\) is replaced by \(S\cup\{x\}\). Explicitly,
\begingroup
\small
\begin{equation}\label{eq:delta-rhs-invariance}
\begin{aligned}
&\rank M[V\setminus S]
-\rank M[V\setminus(S\cup\{y\})]\\
&\quad=\rank M[V\setminus(S\cup\{x\})]
-\rank M[V\setminus(S\cup\{x,y\})].
\end{aligned}
\end{equation}
\endgroup

Since the chosen path is shortest, two vertices on it
are adjacent if and only if they are consecutive along the path.
Thus, $M[S]$ and $M[S\cup\{y\}]$ are adjacency matrices
of paths on $|S|$ and $|S|+1$ vertices, respectively. Elementary row operations and the corresponding column operations reduce each of $M[S]$ and $M[S\cup\{y\}]$ to a block diagonal matrix whose blocks are $\begin{pmatrix}0&1\\1&0\end{pmatrix}$,
with one additional $1\times1$ zero block when its order is odd.
Since these operations preserve rank and each $2\times2$
block has rank two, we obtain
\[
\rank M[S]=2\left\lfloor\frac{|S|}{2}\right\rfloor,
\qquad
\rank M[S\cup\{y\}]
=2\left\lfloor\frac{|S|+1}{2}\right\rfloor.
\]
Thus, adding \(x\) to \(S\) changes the left-hand side from $0$ to $2$ when $|S|$ is even, and from $2$ to $0$
when $|S|$ is odd. However, the right-hand side of
\eqref{eq:delta-root-increment} remains unchanged.
Since both \(S\) and \(S\cup\{x\}\) are subsets of \(X\),
this contradicts \eqref{eq:delta-root-increment}
and proves the claim.

The claim applies to every vertex incident with an edge
joining \(X\) and \(Y\), and also with \(X\) and \(Y\)
interchanged.

Let \(h\in\ker M\) satisfy \(h_y=0\).
Deleting its zero \(y\)-coordinate gives a vector in
\(\ker M[V\setminus\{y\}]\).
By the claim, \(h\) is zero on every component of \(G[X]\)
containing a neighbor of \(y\).

Lemma~\ref{lem:cross-edge-characterization} gives \(h_x=h_y\) for every edge \(xy\) joining \(X\) and \(Y\). Hence, if \(h\) is zero on one component, the claim implies that it is also zero on every component in the opposite part joined to it by an edge. Contracting the components of \(G[X]\) and \(G[Y]\) gives a connected graph. Repeating this argument therefore gives \(h=0\).

We have shown that if \(h \in \ker M\) and \(h_y=0\), then \(h=0\). Thus every nonzero \(h\in\ker M\) satisfies \(h_y=1\). For any two nonzero vectors \(h_1,h_2\in\ker M\), we have \((h_1+h_2)_y=0\), and hence \(h_1+h_2=0\). Thus \(h_1=h_2\), so \(\ker M\) contains at most one nonzero vector. Since $e_y\notin\im M$, the matrix $M$ is singular, and its kernel contains a nonzero vector. Thus $\ker M$ consists of the zero vector and a unique nonzero vector. Over $\mathbb F_2$, this means that $\ker M$ has dimension one. Consequently,
\[
\corank M=1.
\]

Let \(h\) denote the unique nonzero vector in \(\ker M\).
Since $\im M=(\ker M)^\perp$, we have
\[
e_x\notin\im M
\quad\Longleftrightarrow\quad h_x=1,
\]
and
\[
e_x+e_y\in\im M
\quad\Longleftrightarrow\quad h_x=h_y.
\]
Applying Lemma~\ref{lem:cross-edge-characterization}, for $x\in X$
and $y\in Y$ we obtain \(m_{xy}=h_xh_y\). Thus $M[X,Y]=h[X]h[Y]^{\mathsf T}$.
\end{proof}

\begin{lemma}\label{lem:rooted-matrix-uniqueness}
Let \(M\) and \(N\) be symmetric matrices with zero diagonal
over \(\mathbb F_2\), indexed by the same finite set \(V\),
and let \(b\in\mathbb F_2^V\) be a column vector.
Suppose that at least one of the graphs with adjacency matrices
\[
\begin{pmatrix}
0 & b^{\mathsf T}\\
b & M
\end{pmatrix},
\qquad
\begin{pmatrix}
0 & b^{\mathsf T}\\
b & N
\end{pmatrix}
\]
is connected.
If, for every \(S\subseteq V\),
\[
b[S]\in\im M[S]
\quad\Longleftrightarrow\quad
b[S]\in\im N[S],
\]
then \(M=N\).
\end{lemma}

\begin{proof}
We argue by induction on $|V|$. If $|V|\le1$, both $M$ and $N$ are zero matrices, since their diagonal
entries are zero and there are no off-diagonal entries. Hence $M=N$. Now suppose that \(|V|\ge 2\), and assume that the assertion holds for every index set of cardinality less than \(|V|\). Denote the two graphs in the statement by $G_M$ and $G_N$, respectively, and let $r$ be the vertex corresponding to the first row and column of each displayed matrix. Thus both graphs have vertex set $V\cup\{r\}$, where $r\notin V$, and $r$ is adjacent to $v\in V$ if and only if $b_v=1$. Interchanging $M$ and $N$ if necessary, assume that $G_M$ is connected.

Fix $u\in V$ with $b_u=1$, and let $v\in V\setminus\{u\}$ be arbitrary. Then
$b[\{u,v\}]$ is nonzero, and
\[M[\{u,v\}]=
\begin{pmatrix}
0&m_{uv}\\
m_{uv}&0
\end{pmatrix}.
\]
If $m_{uv}=0$, this matrix is zero, so $b[\{u,v\}]\notin\im M[\{u,v\}]$. If $m_{uv}=1$, this matrix is invertible, so $b[\{u,v\}]\in\im M[\{u,v\}]$.
Hence
\[
b[\{u,v\}]\in\im M[\{u,v\}]
\quad\Longleftrightarrow\quad m_{uv}=1.
\]
The same equivalence holds for $N$. The hypothesis therefore gives $m_{uv}=n_{uv}$
for every $v\ne u$. Since $m_{uu}=n_{uu}=0$, the $u$th rows of $M$ and $N$ are identical whenever $b_u=1$.

Thus \(m_{xy}=n_{xy}\) whenever \(b_x=1\) or \(b_y=1\).
It remains to prove the equality when \(b_x=b_y=0\).

Suppose that $M\ne N$, and choose distinct $x,y\in V$ with $m_{xy}\ne n_{xy}$. The preceding observation gives $b_x=b_y=0$. Consider the subgraph of $G_M$ induced by the vertices $v\in V$ with $b_v=0$. Let $S$ be the vertex set of its connected component containing $x$. No vertex of $S$ is adjacent to \(r\), since $b[S]=0$. Moreover, no edge joins $S$ to a vertex $v\in V\setminus S$ with $b_v=0$, since $S$ is a connected component of the induced subgraph on these vertices.
Since $G_M$ is connected, an edge must therefore join $S$ to some vertex $u\in V$ with $b_u=1$.

For every $T\subseteq V$ with $b[T]=0$, order the coordinates of $T\cup\{u\}$ with $u$ first. Then
\[
M[T\cup\{u\}]
=
\begin{pmatrix}
0&M[u,T]\\
M[T,u]&M[T]
\end{pmatrix},
\qquad
b[T\cup\{u\}]
=
\begin{pmatrix}1\\0\end{pmatrix}.
\]
Suppose that $M[T,u]\in\im M[T]$. Here $M[T,u]$ denotes the column vector obtained by
restricting column $u$ of $M$ to the rows indexed by $T$,
and $M[u,T]=M[T,u]^{\mathsf T}$. Choose a column vector $q\in\mathbb F_2^T$ satisfying
\[
M[T]q=M[T,u].
\]
Then $M[u,T]q=q^{\mathsf T}M[T]q=0$, and hence
\[
M[T\cup\{u\}]
\begin{pmatrix}1\\q\end{pmatrix}=0.
\]
Conversely, any kernel vector with $u$-coordinate one
has this form and gives $M[T]q=M[T,u]$.
Thus $M[T,u]\in\im M[T]$ if and only if
$M[T\cup\{u\}]$ has a kernel vector with
$u$-coordinate one.
Since $M[T\cup\{u\}]$ is symmetric,
\[
\im M[T\cup\{u\}]=\bigl(\ker M[T\cup\{u\}]\bigr)^\perp.
\]
Moreover, $b[T\cup\{u\}]$ is the standard basis vector
at $u$, so its inner product with any vector equals
that vector's $u$-coordinate. Thus
\[
b[T\cup\{u\}]\notin\im M[T\cup\{u\}] 
\]
if and only if there is a kernel vector with
$u$-coordinate one.

The same equivalence holds for $N$.
Consequently, the hypothesis implies
\[
M[T,u]\in\im M[T]
\quad\Longleftrightarrow\quad
N[T,u]\in\im N[T]
\qquad(b[T]=0).
\]
Moreover, $M[T,u]=N[T,u]$, since $b_u=1$.

Now $G_M[S\cup\{u\}]$ is connected.
Applying the induction hypothesis to $M[S]$ and $N[S]$,
taking $M[S,u]=N[S,u]$ as the column vector $b$
in the induction hypothesis, gives
\[
M[S]=N[S].
\]
If $y\in S$, then both $x$ and $y$ belong to $S$, so $M[S]=N[S]$ gives $m_{xy}=n_{xy}$, contrary to the choice of $x$ and $y$.

Otherwise, the definition of $S$ gives $m_{xy}=0$, so $n_{xy}=1$. Since \(M[S]=N[S]\) and \(M[u,S]=N[u,S]\), the graphs \(G_M[S\cup\{u\}]\) and \(G_N[S\cup\{u\}]\) are identical. Hence \(G_N[S\cup\{u\}]\) is connected, and the edge $xy$ shows that $G_N[S\cup\{y,u\}]$ is also connected. For every \(T\subseteq S\cup\{y\}\), we have \(b[T]=0\), so the equivalence above applies to every such \(T\). Applying the induction hypothesis to \(M[S\cup\{y\}]\) and \(N[S\cup\{y\}]\), taking \(M[S\cup\{y\},u]=N[S\cup\{y\},u]\) as the common column vector, therefore gives
\[
M[S\cup\{y\}]=N[S\cup\{y\}],
\]
which gives \(m_{xy}=n_{xy}\), contrary to the choice of \(x\) and \(y\). Both applications involve strictly smaller index sets, since $u\notin S\cup\{y\}$. Thus $M=N$.
\end{proof}

\begin{theorem}\label{thm:delta-irreducibility}
Let $G$ be a connected finite simple graph with at least one vertex. Then $P_G^{\langle\delta\rangle}(z;\mathbf u,\mathbf v)$ is irreducible in $K[z^{\pm1},u_x,v_x:x\in V]$ if and only if at least one of the following conditions holds:
\begin{enumerate}
\item The graph obtained from $G$ by complementing the subgraph induced by its even-degree vertices, while leaving all other adjacencies unchanged, is connected.
\item There exist distinct vertices $x,y\in V$ such that the number of their common neighbors is even when both have even degree, or odd when at least one has odd degree.
\end{enumerate}
\end{theorem}

\begin{proof}
We prove the equivalent statement that the polynomial
is reducible if and only if neither condition holds.
Let $M$ be the adjacency matrix of $G$ over $\mathbb F_2$.
Suppose first that $P_G^{\langle\delta\rangle}$ is reducible.
By Lemma~\ref{lem:laurent-factorization-criterion}, there is a nontrivial partition $V=X\sqcup Y$ such that
\begin{equation}\label{eq:delta-final-separation}
f_\delta(S\cup T)+f_\delta(\varnothing)
=
f_\delta(S)+f_\delta(T)
\qquad(S\subseteq X,\ T\subseteq Y).
\end{equation}
Lemma~\ref{lem:delta-rank-one-cut} gives
$\corank M=1$ and
\[
M[X,Y]=h[X]h[Y]^{\mathsf T},
\]
where $h$ is the nonzero vector in $\ker M$.

Taking $S=X$ and $T=Y$ in
\eqref{eq:delta-final-separation} gives
\[
\rank M[X]+\rank M[Y]=\rank M.
\]
Hence
\[
\corank M[X]
+\corank M[Y]=1.
\]
Interchanging $X$ and $Y$ if necessary, assume that
$M[X]$ is nonsingular and
$\corank M[Y]=1$.

Since $G$ is connected and $X,Y$ are nonempty, there exist $x\in X$ and $y\in Y$ with $m_{xy}=1$.
By $M[X,Y]=h[X]h[Y]^{\mathsf T}$, we have \(1=m_{xy}=h_xh_y\), so $h_x=h_y=1$. Hence $h[X]\ne0$ and $h[Y]\ne0$.

Writing $Mh=0$ in block form gives
\[
\begin{pmatrix}
M[X]&M[X,Y]\\
M[Y,X]&M[Y]
\end{pmatrix}
\begin{pmatrix}
h[X]\\
h[Y]
\end{pmatrix}
=0.
\]
The equations indexed by $X$ therefore give
\[
M[X]h[X]+M[X,Y]h[Y]=0.
\]
Substituting $M[X,Y]=h[X]h[Y]^{\mathsf T}$, we obtain
\[
M[X]h[X]+h[X]\bigl(h[Y]^{\mathsf T}h[Y]\bigr)=0.
\]
If $h[Y]^{\mathsf T}h[Y]=0$, the nonsingularity of $M[X]$ would force $h[X]=0$. Therefore,
\[
h[Y]^{\mathsf T}h[Y]=1,
\qquad
M[X]h[X]=h[X].
\]
Since $M[X]$ is symmetric with zero diagonal,
\[
h[X]^{\mathsf T}h[X]=h[X]^{\mathsf T}M[X]h[X]=0.
\]
The equations indexed by $Y$ give
\[
M[Y,X]h[X]+M[Y]h[Y]=0.
\]
Since $M[Y,X]=h[Y]h[X]^{\mathsf T}$ and $h[X]^{\mathsf T}h[X]=0$,
we obtain
\[
M[Y]h[Y]
=
h[Y]\bigl(h[X]^{\mathsf T}h[X]\bigr)=0.
\]
Thus $h[Y]\in\ker M[Y]$.
Since $h[Y]\ne0$ and $\corank M[Y]=1$,
this nonzero vector spans the kernel:
\[
\ker M[Y]=\operatorname{span}\{h[Y]\}.
\]

For $S\subseteq X$ and $T\subseteq Y$, we have
\begin{equation}\label{eq:delta-final-cross-rank}
\begin{aligned}
\rank M[S\cup T]
&=\rank M[S]+\rank M[T]\\
&\quad+\begin{cases}
2,
& h[S]\notin\im M[S]
  \text{ and }h[T]\notin\im M[T],\\
0,
& \text{otherwise}.
\end{cases}
\end{aligned}
\end{equation}
Indeed, with coordinates ordered by $S$ and $T$,
\[
M[S\cup T]
=
\begin{pmatrix}
M[S]&h[S]h[T]^{\mathsf T}\\
h[T]h[S]^{\mathsf T}&M[T]
\end{pmatrix}.
\]
If $h[S]\in\im M[S]$, choose $q$ satisfying
$M[S]q=h[S]$. Then
\[
\begin{pmatrix}
I&0\\
h[T]q^{\mathsf T}&I
\end{pmatrix}
M[S\cup T]
\begin{pmatrix}
I&q h[T]^{\mathsf T}\\
0&I
\end{pmatrix}
=
\begin{pmatrix}
M[S]&0\\
0&M[T]
\end{pmatrix},
\]
since $q^{\mathsf T}h[S]=q^{\mathsf T}M[S]q=0$.
Both multiplying matrices are invertible, so
\[
\rank M[S\cup T]=\rank M[S]+\rank M[T].
\]
The same argument applies when $h[T]\in\im M[T]$.

Suppose now that
$h[S]\notin\im M[S]$ and $h[T]\notin\im M[T]$.
With coordinates ordered by $S$ and $T$, a vector
$\binom{q}{p}$ belongs to $\ker M[S\cup T]$
if and only if
\[
\begin{cases}
M[S]q+h[S]\bigl(h[T]^{\mathsf T}p\bigr)=0,\\
M[T]p+h[T]\bigl(h[S]^{\mathsf T}q\bigr)=0.
\end{cases}
\]
If $h[T]^{\mathsf T}p=1$, the first equation would give
$M[S]q=h[S]$, contrary to $h[S]\notin\im M[S]$.
Hence $h[T]^{\mathsf T}p=0$, and similarly $h[S]^{\mathsf T}q=0$.
The equations therefore reduce to
\[
M[S]q=0,\qquad M[T]p=0,
\qquad h[S]^{\mathsf T}q=0,\qquad h[T]^{\mathsf T}p=0.
\]
Thus
\[
\ker M[S\cup T]
=
\bigl(\ker M[S]\cap h[S]^\perp\bigr)
\oplus
\bigl(\ker M[T]\cap h[T]^\perp\bigr).
\]
Since $M[S]$ is symmetric,
$h[S]\notin\im M[S]=(\ker M[S])^\perp$.
Thus the linear map
\[
\ker M[S]\longrightarrow\mathbb F_2,
\qquad q\longmapsto h[S]^{\mathsf T}q
\]
is nonzero and therefore has rank one.
Its kernel is $\ker M[S]\cap h[S]^\perp$, so
\[
\dim\bigl(\ker M[S]\cap h[S]^\perp\bigr)
=
\corank M[S]-1.
\]
The same argument applies to $M[T]$. Consequently,
\[
\corank M[S\cup T]
=
\corank M[S]
+\corank M[T]-2,
\]
and therefore
\[
\rank M[S\cup T]=\rank M[S]+\rank M[T]+2.
\]

We next determine how membership of the restricted kernel vector in
the column space changes under complementation within each part.
For the part $X$, we prove that
\begin{equation}\label{eq:delta-final-X-membership}
h[S]\in\im M[S]
\quad\Longleftrightarrow\quad
h[X\setminus S]\in\im M[X\setminus S]
\qquad(S\subseteq X).
\end{equation}
Since $h[Y]\in\ker M[Y]$ and $h[Y]^{\mathsf T}h[Y]=1$, \(h[Y]\notin\im M[Y]\). Also, \(f_\delta(Y)=f_\delta(\varnothing)\), so \eqref{eq:delta-final-separation} gives \(f_\delta(S\cup Y)=f_\delta(S)\). By \eqref{eq:delta-final-cross-rank}, both sides have the common terms
\[
\rank M[S]+\rank M[X\setminus S]+\rank M[Y].
\]
The left-hand side has an additional $2$ precisely when $h[S]\notin\im M[S]$, whereas the right-hand side has an additional $2$ precisely when $h[X\setminus S]\notin\im M[X\setminus S]$. This proves \eqref{eq:delta-final-X-membership}.

We also have
\begin{equation}\label{eq:delta-final-Y-membership}
h[T]\in\im M[T]
\quad\Longleftrightarrow\quad
h[Y\setminus T]\notin\im M[Y\setminus T]
\qquad(T\subseteq Y).
\end{equation}
Choose $x\in X$ with $h_x=1$. Since $M[\{x\}]=0$, \eqref{eq:delta-final-X-membership} gives
\[
h[X\setminus\{x\}]
\notin\im M[X\setminus\{x\}].
\]
The nonsingularity of $M[X]$ and Lemma~\ref{lem:alternating-deletion} imply
\[
\rank M[X\setminus\{x\}]=\rank M[X]-2.
\]
Moreover, if $e_x=Mq$ for some $q$, then
\[
h_x=h^{\mathsf T}e_x=h^{\mathsf T}Mq=(Mh)^{\mathsf T}q=0,
\]
contrary to $h_x=1$. Hence $e_x\notin\im M$.
Lemma~\ref{lem:alternating-deletion} therefore gives \(f_\delta(\{x\})=f_\delta(\varnothing)\). Consequently, \eqref{eq:delta-final-separation} yields \(f_\delta(\{x\}\cup T)=f_\delta(T)\).

Since $M[X]$ is nonsingular, we have $h[X]\in\im M[X]$. \eqref{eq:delta-final-cross-rank} gives
\[
f_\delta(T)
=
\rank M[X]+\rank M[T]+\rank M[Y\setminus T].
\]
Since $h_x=1$ and $h[X\setminus\{x\}]\notin\im M[X\setminus\{x\}]$, the definition of \(f_{\delta}\) gives
\[
f_\delta(\{x\}\cup T)
=
\rank M[\{x\}\cup T]
+\rank M[(X\setminus\{x\})\cup(Y\setminus T)].
\]
Using \eqref{eq:delta-final-cross-rank} and $\rank M[X\setminus\{x\}]=\rank M[X]-2$, we obtain
\[
\begin{aligned}
f_\delta(\{x\}\cup T)
={}&f_\delta(T)-2\\
&+\begin{cases}
2,&h[T]\notin\im M[T],\\
0,&h[T]\in\im M[T],
\end{cases}\\
&+\begin{cases}
2,&h[Y\setminus T]\notin\im M[Y\setminus T],\\
0,&h[Y\setminus T]\in\im M[Y\setminus T].
\end{cases}
\end{aligned}
\]
The equality $f_\delta(\{x\}\cup T)=f_\delta(T)$ therefore implies that exactly one of
\[
h[T]\notin\im M[T],
\qquad
h[Y\setminus T]\notin\im M[Y\setminus T]
\]
holds. This proves \eqref{eq:delta-final-Y-membership}.

It follows that, for every $A\subseteq V$,
\begin{equation}\label{eq:delta-final-complement-membership}
h[A]\in\im M[A]
\quad\Longleftrightarrow\quad
h[A^c]\notin\im M[A^c].
\end{equation}
To see this, write $A=S\cup T$, where $S\subseteq X$
and $T\subseteq Y$. Then
\[
M[A]+h[A]h[A]^{\mathsf T}
=
\begin{pmatrix}
M[S]+h[S]h[S]^{\mathsf T}&0\\
0&M[T]+h[T]h[T]^{\mathsf T}
\end{pmatrix}.
\]
By \eqref{eq:delta-final-cross-rank} and
Lemma~\ref{lem:rank-one-update}, the rank of this matrix
equals $\rank M[A]$ if both
\[
h[S]\in\im M[S],
\qquad
h[T]\in\im M[T]
\]
hold or neither holds, and exceeds $\rank M[A]$ by one
if exactly one holds.
Applying Lemma~\ref{lem:rank-one-update} to $M[A]$,
we conclude that $h[A]\in\im M[A]$ precisely when
both conditions hold or neither holds.

Under complementation,
\eqref{eq:delta-final-X-membership} leaves the first
condition unchanged, whereas
\eqref{eq:delta-final-Y-membership} reverses the second.
Thus exactly one condition holds for $A^c$ if both
hold or neither holds for $A$, and conversely.
This proves \eqref{eq:delta-final-complement-membership}.

Since $h\in\ker M$ and $h^{\mathsf T}h=1$, we have
$h\notin\im M$. Lemma~\ref{lem:rank-one-update}
therefore shows that $M+hh^{\mathsf T}$ is nonsingular.
Define
\[
N=(M+hh^{\mathsf T})^{-1}+hh^{\mathsf T}.
\]
We will apply Lemma~\ref{lem:rooted-matrix-uniqueness} to show that $N=M$.
This will turn the identity $MN=I+hh^{\mathsf T}$ established below into
the required identity for $M^2$.
Since $(M+hh^{\mathsf T})h=h$, we have $(M+hh^{\mathsf T})^{-1}h=h$. Thus
\[
Nh=(M+hh^{\mathsf T})^{-1}h+h(h^{\mathsf T}h)=h+h=0.
\]
Using $Mh=0$, we obtain
\[
\begin{aligned}
MN
&=M(M+hh^{\mathsf T})^{-1}\\
&=\bigl((M+hh^{\mathsf T})+hh^{\mathsf T}\bigr)(M+hh^{\mathsf T})^{-1}\\
&=I+hh^{\mathsf T}.
\end{aligned}
\]
Since $M$ and $N$ are symmetric, transposing this identity also gives $NM=I+hh^{\mathsf T}$. The identities \(NM=I+hh^{\mathsf T}\) and \(h^{\mathsf T}N=0\) give $N=NMN$. The matrix $N$ is symmetric, and \(q^{\mathsf T}Nq=(Nq)^{\mathsf T}M(Nq)=0\) for every column vector $q$. Thus $N$ has zero diagonal.

We claim that
\[
h[A]\notin\im M[A]
\quad\Longleftrightarrow\quad
h[A^c]\in\im N[A^c]
\qquad(A\subseteq V).
\]
Suppose first that $h[A]\notin\im M[A]$. Since $M[A]$ is symmetric, \(\im M[A]=(\ker M[A])^\perp\). As $h[A]\notin\im M[A]$, choose $q[A]\in\ker M[A]$ with \(h[A]^{\mathsf T}q[A]=1\), and extend it to $q\in\mathbb F_2^V$ by setting $q[A^c]=0$. Thus
\[
q[A^c]=0,\qquad M[A]q[A]=0,\qquad h^{\mathsf T}q=1.
\]
Then $(Mq)[A]=0$ and \(NMq=(I+hh^{\mathsf T})q=q+h\). Restricting to $A^c$ gives \(N[A^c](Mq)[A^c]=h[A^c]\), so $h[A^c]\in\im N[A^c]$.

Conversely, choose $q\in\mathbb F_2^V$ satisfying \(q[A]=0, N[A^c]q[A^c]=h[A^c]\). Since $N$ is symmetric with zero diagonal, \(h^{\mathsf T}q=q^{\mathsf T}Nq=0\). Consequently,
\[
(Nq+h)[A^c]=0,
\qquad
M(Nq+h)=(I+hh^{\mathsf T})q=q.
\]
Thus \(M[A](Nq+h)[A]=q[A]=0\). Moreover, since $Nh=0$ and $h^{\mathsf T}h=1$,
\[
h[A]^{\mathsf T}(Nq+h)[A]
=h^{\mathsf T}(Nq+h)=1.
\]
Hence $(Nq+h)[A]\in\ker M[A]$. Together with $h[A]^{\mathsf T}(Nq+h)[A]=1$, symmetry gives $h[A]\notin\im M[A]$. Combining this equivalence with
\eqref{eq:delta-final-complement-membership} gives
\[
h[A]\in\im M[A]
\quad\Longleftrightarrow\quad
h[A]\in\im N[A]
\qquad(A\subseteq V).
\]
Since $G$ is connected and $h\ne0$, the graph with
adjacency matrix
\[
\begin{pmatrix}
0&h^{\mathsf T}\\
h&M
\end{pmatrix}
\]
is connected.
Lemma~\ref{lem:rooted-matrix-uniqueness} gives $N=M$.
Hence
\begin{equation}\label{eq:delta-final-square}
M^2=I+hh^{\mathsf T}.
\end{equation}
Comparison of diagonal entries gives
\[
\deg_G(x)\equiv1+h_x\pmod2
\qquad(x\in V).
\]
Thus $h$ is the indicator vector of the even-degree
vertices of $G$.
For distinct $x,y\in V$, comparison of off-diagonal
entries gives
\[
|N_G(x)\cap N_G(y)|\equiv h_xh_y\pmod2,
\]
which shows that condition~\textup{(2)} does not hold.

By $M[X,Y]=h[X]h[Y]^{\mathsf T}$, we have $m_{xy}=h_xh_y$ for every $x\in X$ and $y\in Y$. Since $h$ is the indicator vector of the even-degree vertices, two vertices in opposite parts are adjacent if and only if both have even degree. Complementing the subgraph induced by the even-degree vertices therefore removes every edge between $X$ and $Y$ and adds none between them. Since $X$ and $Y$ are nonempty, the resulting graph is disconnected. Thus condition~\textup{(1)} does not hold either.

Conversely, suppose that neither condition holds.
Let $h\in\mathbb F_2^V$ be the column vector with
$h_x=1$ when $\deg_G(x)$ is even and $h_x=0$ otherwise.
The diagonal entries of $M^2$ are the degree parities,
and the failure of condition~\textup{(2)} determines its off-diagonal
entries. Hence
\[
M^2=I+hh^{\mathsf T}.
\]
Also, $h\ne0$, since otherwise the graph in
condition~\textup{(1)} would coincide with the connected graph $G$,
contrary to the failure of that condition.

By associativity, $MM^2=M^2M$.
Substituting $M^2=I+hh^{\mathsf T}$ gives
\[
M+(Mh)h^{\mathsf T}=M+h(h^{\mathsf T}M).
\]
Thus
\[
(Mh)h^{\mathsf T}=h(h^{\mathsf T}M)=h(Mh)^{\mathsf T}.
\]
Since $h\ne0$, choose $x\in V$ with $h_x=1$. The $x$th columns of $(Mh)h^{\mathsf T}=h(Mh)^{\mathsf T}$ give
\[
(Mh)h_x=h(Mh)_x.
\]
Thus $Mh=(Mh)_x h$.
Since $M$ is symmetric with zero diagonal over
$\mathbb F_2$, we have $q^{\mathsf T}Mq=0$ for every column
vector $q$. Taking $q=Me_x$ gives
\[
(M^3)_{xx}
=e_x^{\mathsf T}M^3e_x
=(Me_x)^{\mathsf T}M(Me_x)=0.
\]
Since $M^3=M+(Mh)h^{\mathsf T}$, we obtain $(Mh)_x=0$,
and hence $Mh=0$.

Applying $M^2=I+hh^{\mathsf T}$ to $h$ gives
$0=h+h(h^{\mathsf T}h)$.
Since $h\ne0$, it follows that $h^{\mathsf T}h=1$.
Consequently,
\[
(M+hh^{\mathsf T})^2=I.
\]
Since $Mh=0$ and $h\ne0$, we have $\rank M\le |V|-1$.
On the other hand, $(M+hh^{\mathsf T})^2=I$ and $\rank(hh^{\mathsf T})=1$
give
\[
|V|=\rank(M+hh^{\mathsf T})\le\rank M+1.
\]
Hence $\rank M=|V|-1$, so $\corank M=1$.
Since $M$ is symmetric with zero diagonal, its rank is
even. Therefore $|V|=\rank M+1$ is odd.

By \cite[Proposition~2]{Brijder2013Nullity},
complementary principal submatrices of an invertible
matrix and its inverse have equal corank.
Since $(M+hh^{\mathsf T})^{-1}=M+hh^{\mathsf T}$, we obtain
\[
\operatorname{corank}\bigl((M+hh^{\mathsf T})[S]\bigr)
=
\operatorname{corank}\bigl((M+hh^{\mathsf T})[S^c]\bigr)
\qquad (S\subseteq V).
\]
If $M[A]q=h[A]$, then $h[A]^{\mathsf T}q=q^{\mathsf T}M[A]q=0$, since $M[A]$ is symmetric with zero diagonal. Thus Lemma~\ref{lem:rank-one-update} shows that the rank cannot decrease and gives
\[
\begin{gathered}
\corank M[A]\in
\left\{
\corank\bigl((M+hh^{\mathsf T})[A]\bigr),
\corank\bigl((M+hh^{\mathsf T})[A]\bigr)+1
\right\}\\
(A\subseteq V).
\end{gathered}
\]
Both $M[S]$ and $M[S^c]$ have even rank.
Since $|V|$ is odd, their coranks have opposite parity.
The preceding two relations therefore imply that
one corank equals
$\corank((M+hh^{\mathsf T})[S])$
and the other exceeds it by one. Hence
\[
\corank M[S]
+\corank M[S^c]
=
1+2\corank\bigl((M+hh^{\mathsf T})[S]\bigr).
\]
Therefore,
\begin{equation}\label{eq:delta-final-exponent}
f_\delta(S)
=
|V|-1
-2\corank\bigl((M+hh^{\mathsf T})[S]\bigr).
\end{equation}

The off-diagonal entries of $M+hh^{\mathsf T}$ are the adjacency
entries of the graph in condition~\textup{(1)}.
Let $X$ be the vertex set of one of its connected
components, and put $Y=V\setminus X$.
Since condition~\textup{(1)} does not hold, this graph is
disconnected, so both sets are nonempty.
Since $(M+hh^{\mathsf T})[X,Y]=0$,
\[
\corank\bigl((M+hh^{\mathsf T})[S\cup T]\bigr)
=
\corank\bigl((M+hh^{\mathsf T})[S]\bigr)
+
\corank\bigl((M+hh^{\mathsf T})[T]\bigr)
\]
for $S\subseteq X$ and $T\subseteq Y$.
Together with \eqref{eq:delta-final-exponent}, this gives
\[
f_\delta(S\cup T)+f_\delta(\varnothing)
=
f_\delta(S)+f_\delta(T).
\]
Lemma~\ref{lem:laurent-factorization-criterion}
therefore implies that $P_G^{\langle\delta\rangle}$
is reducible.
\end{proof}

\begin{corollary}\label{cor:delta-irreducible-factorization}
Let $G$ be a connected finite simple graph for which neither
condition in Theorem~\ref{thm:delta-irreducibility} holds.
Let $M$ be its adjacency matrix over $\mathbb F_2$,
and let $h\in\mathbb F_2^V$ be the column vector with
$h_x=1$ when $\deg_G(x)$ is even and $h_x=0$ otherwise.
Let $V_1,\ldots,V_t$ be the vertex sets of the connected
components of the graph in condition~\textup{(1)}.
Then
\[
P_G^{\langle\delta\rangle}(z;\mathbf u,\mathbf v)
=
z^{|V|-1}
\prod_{i=1}^{t}
\left(
\sum_{S\subseteq V_i}
u_Sv_{V_i\setminus S}\,
z^{-2\corank\bigl((M+hh^{\mathsf T})[S]\bigr)}
\right),
\]
and every factor in parentheses is irreducible in
$K[z^{\pm1},\mathbf u,\mathbf v]$.
\end{corollary}

\begin{proof}
Put $B=M+hh^{\mathsf T}$. Its block diagonal form gives
\[
\corank B[A]
=\sum_{i=1}^{t}\corank B[A\cap V_i].
\]
Together with \eqref{eq:delta-final-exponent}, this yields the stated factorization by summing
independently over the subsets $A\cap V_i$.

Fix $i$ and suppose that its factor is reducible. Factors cannot involve vertex variables outside $V_i$, by additivity of degree. Apply Lemma~\ref{lem:laurent-factorization-criterion}
after multiplying by the unit $z^{2|V_i|}$. The exponent function is
$g(S)=2|V_i|-2\corank B[S]\geq0$. Thus there is a nontrivial partition $V_i=X\sqcup Y$
such that
\[
g(S\cup T)+g(\varnothing)=g(S)+g(T)
\qquad(S\subseteq X,\ T\subseteq Y).
\] Substituting $g$ and cancelling $4|V_i|$ and then $-2$
gives
\[
\corank B[S\cup T]
=\corank B[S]+\corank B[T]
\qquad(S\subseteq X,\ T\subseteq Y).
\]

The connected component with vertex set $V_i$ has an edge $xy$
across this partition.
Since $B_{xx}=h_x$, $B_{yy}=h_y$, and $B_{xy}=1$,
taking $S=\{x\}$ and $T=\{y\}$ gives
\[
2-\operatorname{rank}
\begin{pmatrix}
h_x&1\\
1&h_y
\end{pmatrix}
=(1-h_x)+(1-h_y).
\]
Hence
\[
\operatorname{rank}
\begin{pmatrix}
h_x&1\\
1&h_y
\end{pmatrix}
=h_x+h_y,
\]
where $h_x,h_y$ on the right are viewed as integers.
If $h_x=h_y=1$, the matrix has two equal nonzero rows,
so its rank is one rather than two.
Otherwise, its determinant over $\mathbb F_2$ is
$h_xh_y-1=1$, so its rank is two,
whereas $h_x+h_y\leq1$.
This contradiction proves irreducibility.
\end{proof}

\begin{example}
Let $G=K_{1,3}$ have vertex set $V=\{1,2,3,4\}$ and edge
set $\{12,13,14\}$. From the definition,
\[
P_{K_{1,3}}^{\langle\delta\rangle}(z;\mathbf u,\mathbf v)
=z^2\prod_{x=1}^{4}(u_x+v_x)
+(1-z^2)\bigl(u_1v_2v_3v_4+v_1u_2u_3u_4\bigr).
\]
All vertices of $K_{1,3}$ have odd degree. Thus the transformation in
condition~\textup{(1)} of Theorem~\ref{thm:delta-irreducibility} leaves the
graph unchanged, and the resulting graph is connected. That theorem
therefore shows that this polynomial is irreducible in
$K[z^{\pm1},u_x,v_x:x\in V]$.

Let $G$ have vertex set
$V=\{1,2,3,4,5\}$ and edge set
\[
\{12,13,14,15,25,34\}.
\]
Every vertex has even degree, so the transformation in
condition~\textup{(1)} of Theorem~\ref{thm:delta-irreducibility} complements
the entire graph. The resulting graph has edge set
$\{23,24,35,45\}$ and connected components with vertex sets $\{1\}$ and
$\{2,3,4,5\}$. Every pair of distinct vertices of $G$ has exactly one common
neighbor: for adjacent vertices it is the third vertex of their triangle,
and for nonadjacent vertices it is vertex $1$. Thus neither condition in
Theorem~\ref{thm:delta-irreducibility} holds.
The definition gives
\[
\begin{aligned}
P_G^{\langle\delta\rangle}(z;\mathbf u,\mathbf v)
={}&z^2(u_1+v_1)\Bigl(
z^2\prod_{x=2}^{5}(u_x+v_x)\\
&\quad +(1-z^2)(
u_2u_3v_4v_5+u_2u_4v_3v_5
+u_3u_5v_2v_4+u_4u_5v_2v_3)
\Bigr).
\end{aligned}
\]
This is a nontrivial factorization in the Laurent polynomial
ring. The factor indexed by $\{2,3,4,5\}$ is irreducible by
Corollary~\ref{cor:delta-irreducible-factorization}; the factor $u_1+v_1$
is linear and hence irreducible.
\end{example}

\subsection{Irreducibility in the ordinary polynomial ring}
The preceding results admit the following unified formulation.
Although the five exponent functions have different matrix
expressions, irreducibility in the ordinary polynomial ring is
governed in every case by the same two conditions: connectedness
of the underlying graph and the existence of a state with
zero \(z\)-exponent.

\begin{theorem}\label{thm:unified-ordinary-irreducibility}
Let $G$ be a finite simple graph with nonempty vertex
set $V$, and let the operation be
\[
\bullet\in\{\delta,\tau,\delta\tau,\tau\delta,\delta\tau\delta\}.
\]
Then $P_G^{\langle\bullet\rangle}(z;\mathbf u,\mathbf v)$
is irreducible in $K[z,u_x,v_x:x\in V]$ if and only if
$G$ is connected and
\[
\min_{A\subseteq V}f_\bullet(A)=0.
\]
\end{theorem}

\begin{proof}
By Corollary~\ref{cor:ordinary-factorization-criterion}
and Lemma~\ref{lem:laurent-factorization-criterion},
$P_G^{\langle\bullet\rangle}$ is irreducible in the
ordinary polynomial ring if and only if
\[
\min_{A\subseteq V}f_\bullet(A)=0
\]
and it is irreducible in
$K[z^{\pm1},u_x,v_x:x\in V]$.

For
$\bullet\in\{\tau,\delta\tau,\tau\delta,\delta\tau\delta\}$,
Theorems~\ref{thm:tau-irreducibility-laurent},
\ref{thm:delta-tau-irreducibility-laurent},
\ref{thm:tau-delta-irreducibility}, and
\ref{thm:delta-tau-delta-laurent-irreducibility},
respectively, show that Laurent irreducibility is
equivalent to the connectedness of $G$.
Thus the assertion follows in these four cases.

It remains to consider $\bullet=\delta$.
Since
\[
f_\delta(A)=\rank M[A]+\rank M[A^c],
\]
we have $f_\delta(A)=0$ precisely when both $A$ and
$A^c$ are independent sets. Consequently,
\[
\min_{A\subseteq V}f_\delta(A)=0
\quad\Longleftrightarrow\quad
G\text{ is bipartite}.
\]

If $G$ is disconnected, its connected components give
a nontrivial factorization of $P_G^{\langle\delta\rangle}$
into nonunits, so it is reducible in the Laurent
polynomial ring.

Now suppose that $G$ is connected and bipartite.
If condition~\textup{(2)} of
Theorem~\ref{thm:delta-irreducibility} holds, that theorem
already gives Laurent irreducibility.
Otherwise, for all distinct $x,y\in V$, the number of
their common neighbors is odd if and only if both
vertices have even degree.

Vertices in opposite bipartition classes have no common
neighbors, so they cannot both have even degree.
Hence all even-degree vertices lie in one bipartition
class, and the subgraph they induce is edgeless.
Complementing this subgraph only adds edges to $G$.
The resulting graph is therefore connected, so
condition~\textup{(1)} of
Theorem~\ref{thm:delta-irreducibility} holds.
That theorem again gives Laurent irreducibility.

Thus, when $\min_{A\subseteq V}f_\delta(A)=0$,
Laurent irreducibility is equivalent to the connectedness
of $G$. Combining this with the ordinary polynomial
ring criterion proves the assertion for $\bullet=\delta$
as well.
\end{proof}

\begin{corollary}
Let $G$ be a finite simple graph with at least one vertex,
and let $\bullet\in\{\tau,\delta\tau\}$.
Then $P_G^{\langle\bullet\rangle}(z;\mathbf u,\mathbf v)$
is irreducible in
$K[z,u_x,v_x:x\in V(G)]$
if and only if $G$ has exactly one vertex.
\end{corollary}

\begin{proof}
For $\bullet=\tau$, the equality
$\operatorname{rank}(M+I_A)=0$ implies $M+I_A=0$.
Since $M$ has zero diagonal, this forces $A=\varnothing$
and $M=0$.

For $\bullet=\delta\tau$, a zero exponent requires
\[
\operatorname{rank}(M[A]+I_{|A|})
=\operatorname{rank}M[A^c]=0.
\]
The diagonal entries of $M[A]+I_{|A|}$ are all one,
so $A=\varnothing$, and then $M=0$.
Conversely, if $M=0$, the empty subset has exponent zero
in both cases.

Thus, in either case, a zero exponent exists precisely
when $G$ is edgeless. The conclusion now follows
from Theorem~\ref{thm:unified-ordinary-irreducibility}.
\end{proof}

\bibliographystyle{elsarticle-num}
\bibliography{references}

\end{document}